\pdfoutput=1
\RequirePackage{ifpdf}
\ifpdf 
\documentclass[pdftex]{sigma}
\else
\documentclass{sigma}
\fi

\numberwithin{equation}{section}

\usepackage{multirow}

\newtheorem{Theorem}{Theorem}[section]

\newtheorem{Lemma}[Theorem]{Lemma}
\newtheorem{Proposition}[Theorem]{Proposition}
 { \theoremstyle{definition}

 }

\newcommand{\dtime}{\frac{{\rm d}}{{\rm d} t}} 
\newcommand{\dvar}[1]{\partial_{#1}} 

\newcommand{\tr}{\operatorname{tr}} 
\newcommand{\grad}{\operatorname{grad}} 
\newcommand{\diag}{\operatorname{diag}} 

\begin{document}
\allowdisplaybreaks

\newcommand{\arXivNumber}{1803.07733}

\renewcommand{\PaperNumber}{027}

\FirstPageHeading

\ShortArticleName{Bach Flow on Homogeneous Products}

\ArticleName{Bach Flow on Homogeneous Products}

\Author{Dylan HELLIWELL}

\AuthorNameForHeading{D.~Helliwell}

\Address{Department of Mathematics, Seattle University, 901 12th Ave, Seattle, WA 98122, USA}
\Email{\href{mailto:helliwed@seattleu.edu}{helliwed@seattleu.edu}}

\ArticleDates{Received September 03, 2019, in final form March 29, 2020; Published online April 11, 2020}

\Abstract{Qualitative behavior of Bach flow is established on compact four-dimensional locally homogeneous product manifolds. This is achieved by lifting to the homogeneous universal cover and, in most cases, capitalizing on the resultant group structure. The resulting system of ordinary differential equations is carefully analyzed on a case-by-case basis, with explicit solutions found in some cases. Limiting behavior of the metric and the curvature are determined in all cases. The behavior on quotients of~$\mathbb{R} \times \mathbb{S}^3$ proves to be the most challenging and interesting.}


\Classification{53C44; 53C30; 34C40}

\vspace{-2mm}

\section{Introduction}

In four dimensions, Bach flow is a solution to
\begin{gather*}
\dvar{t} g = B + \frac{1}{12} \Delta S g, \qquad g(0) = h,
\end{gather*}
where $B$ is the Bach tensor and $S$ is scalar curvature for the metric~$g$. This serves as a concrete motivating example of a higher-order intrinsic curvature flow. Such flows, including flow by the ambient obstruction tensor and flow by the gradient of the total curvature energy functional, have been of interest recently. See for example \cite{BahuaudHelliwell1,BahuaudHelliwell2}, and with related work found in \cite{Bour,Ho,Lopez, Streets}.

Analyzing Bach flow on an arbitrary locally homogeneous 4-manifold is a challenging endeavor, so our goal here is to understand Bach flow on a more restricted family that is more tractable. Specifically, we study Bach flow on $(M, g)$ where $M = \mathbb{S}^1 \times N$, $(N, \tilde{g})$ is a closed locally homogeneous three-dimensional Riemannian manifold, and $g = g_{\mathbb{S}^1} + \tilde{g}$ is the product metric. By lifting to the universal cover $\widehat{M}$ of $M$, this analysis reduces to analysis of Bach flow on one of nine simply connected homogeneous spaces.

The specific details for each of the nine cases can be found in Sections \ref{13flowsec} and \ref{22flowsec}. As a summary, we find:
\begin{itemize}\itemsep=0pt
\item if $\widehat{N} = \mathbb{R}^3$ or $\mathbb{H}^3$, Bach flow is static;
\item if $\widehat{N} = {\rm Nil}$, $\widehat{{\rm SL}}(2,\mathbb{R})$, $\mathbb{R} \times \mathbb{S}^2$, or $\mathbb{R} \times \mathbb{H}^2$, Bach flow collapses to a flat surface;
\item if $\widehat{N} = {\rm Solv}$, Bach flow collapses to a curve;
\item if $\widehat{N} = {\rm E}(2)$, Bach flow converges to a flat four-dimensional manifold;
\item if \looseness=-1 $\widehat{N} = \mathbb{S}^3$, Bach flow can collapse to a flat three-dimensional manifold, collapse to a~flat surface, or converge to a curved four-dimensional manifold, depending on the initial conditions.
\end{itemize}
In this paper, $\mathbb{S}^n$ is the $n$-dimensional sphere, $\mathbb{H}^n$ is $n$-dimensional hyperbolic space, ${\rm Nil}$ is the Heisenberg group consisting of $3 \times 3$ upper triangular matrices with 1's on the diagonal, ${\rm Solv}$ is the Poincar\'{e} group for 2-D Minkowski space $\mathbb{R}^2 \rtimes {\rm O}(1,1)$, ${\rm E}(2)$ is the group of Euclidean transformations of the plane, and $\widehat{{\rm SL}}(2, \mathbb{R})$ is the universal cover of ${\rm SL}(2, \mathbb{R})$.

 The method here is similar to that of \cite{IJ} and \cite{IJL}, where the qualitative behavior of volume-normalized Ricci flow on locally homogeneous three- and four-dimensional manifolds was determined. See~\cite{GlickPayne} for an alternative approach to analyzing Ricci flow on homogeneous three-dimensional manifolds, and see \cite{KM} for analysis of the quasi-convergence of locally homogeneous manifolds under Ricci flow. Other geometric flows have also been analyzed on locally homogeneous spaces. See for example \cite{KST} for analysis of Cotton flow, \cite{CNSC} and \cite{CSC2} for analysis of cross curvature flow, \cite{CSC} for analysis of backward Ricci flow, and~\cite{GGI} for analysis of second-order renormalization group flow.

Here, for most spaces, the analysis is a bit more challenging than for Ricci flow, since the polynomials in the systems determined by Ricci flow are third order while those for Bach flow are seventh order. These higher order expressions are more difficult to analyze for the purposes of qualitative analysis.

Determining the behavior of Bach flow on model spaces has so far been limited to flow on locally homogeneous $2 \times 2$ products by \cite{DasKar} and $\mathbb{S}^1 \times {\rm Solv}$ by~\cite{Ho}. Additionally, in~\cite{Streets}, flow by the gradient of the total curvature energy functional was analyzed on two specific four-dimensional homogeneous spaces: $\mathbb{S}^2 \times \mathbb{H}^2$ and $\mathbb{R} \times \mathbb{S}^3$. Bach flow is related to this flow, and comparing and contrasting the qualitative behavior of these flows helps to understand this relationship. On $\mathbb{S}^2 \times \mathbb{H}^2$ the equations determined by the two flows are essentially the same. On $\widehat{M} = \mathbb{R} \times \mathbb{S}^3$, only round metrics on $\mathbb{S}^3$ were considered in~\cite{Streets}, and on compact quotients, the resulting product metric was found to collapse to a three-dimensional space, with the $\mathbb{S}^1$ slice shrinking. Here, we find that Bach flow is static in this case.

The general approach to understanding Bach flow on the spaces of interest is similar to that found in~\cite{IJ}. The universal cover $\widehat{N}$ is either a Lie group or it is not. In the case where $\widehat{N}$ is a Lie group, the set of homogeneous metrics can be identified with the set of left-invariant metrics on $\widehat{N}$, which in turn are identified with the set of inner products on the tangent space at the identity. Curvature can then be expressed in terms of the structure constants and the inner product. The Lie groups of interest have the property that a basis can be found where the inner product is diagonal and the structure constants can be written in a convenient form. As was true for the Ricci tensor in~\cite{IJ}, we find here that the Bach tensor in this setting is diagonal and so Bach flow preserves this structure. The resulting system is analyzed, with explicit solutions found in some cases, and limiting behavior is determined. If $\widehat{N}$ is not a Lie group, the analysis proves to be somewhat simpler, owing to the fact that there are fewer homogeneous metrics on these spaces. The resulting systems can all be solved explicitly.

This paper is organized as follows: In Section~\ref{Bachbackgroundsec}, the Bach tensor and Bach flow are discussed. Additionally, formulas for the Bach tensor on products are provided. In Section~\ref{spacebackgroundsec}, details surrounding the locally homogeneous spaces and Lie groups of interest are provided, including curvature formulas in terms of structure constants. In Section~\ref{odesec}, useful results about ordinary differential equations are provided. Section~\ref{13flowsec} is devoted to the derivation and analysis of Bach flow on locally homogeneous $1 \times 3$ products and in Section \ref{22flowsec}, Bach flow is analyzed on locally homogeneous $2 \times 2$ products. Finally, in Section~\ref{riccompsec}, the qualitative results for Bach flow are compared and contrasted with those of Ricci flow.

\section{The Bach tensor and Bach flow} \label{Bachbackgroundsec}

On a four-dimensional Riemannian manifold $(M, g)$, the Bach tensor $B$ is given by
\[
B_{jk} = g^{lq} P_{jk;lq} - g^{lq} P_{jl;kq} + P^{il} W_{ijkl},
\]
where $P$ is the Schouten tensor, which, for an $n$-dimensional manifold is defined as
\[
P = \frac{1}{n-2}\left({\rm Ric}- \frac{S}{2(n-1)}g \right)
\]
and $W$ is the Weyl tensor. Throughout this paper, curvature and index conventions follow those found in \cite{Jack3}. The Bach tensor is a symmetric, trace free, divergence free tensor that is fourth order in the metric and is conformally invariant: if $\bar{g} = \rho^2 g$, then $\bar{B} = \rho^{-2} B$. It can be realized as $-\frac{1}{4} \grad(\mathcal{W})$ where $\mathcal{W}$ is the Weyl energy functional:
\[
\mathcal{W} = \int_M |W|^2 {\rm d} \mu
\]
with $|W|^2 = g^{ip} g^{jq} g^{kr} g^{ls} W_{ijkl} W_{pqrs}$.

In \cite{BahuaudHelliwell1} and \cite{BahuaudHelliwell2} short time existence and uniqueness are established for solutions to the geometric flow
\begin{equation*}
\dvar{t} g = B + \frac{1}{12} \Delta S g, \qquad g(0) = h.
\end{equation*}
Here, and throughout, $\Delta = g^{ij} \nabla_i \nabla_j$. The positive multiple of $\Delta S g$ is included to ensure that the resulting flow is well posed. The fraction $\frac{1}{12}$ could be replaced by any positive constant $\alpha$. In general, if $\alpha = 0$, the analysis in~\cite{BahuaudHelliwell1} and~\cite{BahuaudHelliwell2} no longer applies and, as far as the author is aware, it is not known if the case $\alpha = 0$ is well posed.

On the other hand if, for a solution to the flow above, the scalar curvature $S(t)$ is constant on $M$, then the flow reduces to
\begin{gather} \label{bachfloweqn}
\dvar{t} g = B, \qquad
g(0) = h.
\end{gather}
In this paper, the flow is analyzed on locally homogeneous product manifolds and local homogeneity ensures that the scalar curvature is constant on~$M$.

One useful consequence of the fact that the Bach tensor is trace free is that the volume form is constant along the flow:
\[
\dvar{t}{\rm d} \mu_g = 0.
\]

If fixed coordinates are chosen, this is equivalent to saying $\det g(t) = \det h$ is constant in time.

\subsection{Bach tensor on products}

In general, the Bach flow equations lead to a complicated nonlinear system. Making use of the product structure significantly simplifies the resulting equations.

Let $\big(N^{(1)}, \tilde{g}^{(1)}\big)$ and $\big(N^{(2)}, \tilde{g}^{(2)}\big)$ be Riemannian manifolds. Let $M = N^{(1)} \times N^{(2)}$. The product metric $g$ on $M$ is
\begin{gather*}
g = g^{(1)} + g^{(2)},
\end{gather*}
where $g^{(i)}= \pi_i^*\big(\tilde{g}^{(i)}\big)$ are the pullbacks of the component metrics by the natural projections. Greek indices ($\alpha$, $\beta$, $\gamma$, etc.) will be used for $N^{(1)}$, and lower case roman indices ($i$, $j$, $k$, etc.) will be used for $N^{(2)}$. In the case where $N^{(1)}$ is one-dimensional, the subscript $0$ will be used. Abusing notation slightly, the tildes used above will be dropped. To clarify when dealing with an object on $N^{(1)}$ or $N^{(2)}$ (as opposed to $M$) a parenthetical superscript will be used to indicate the component.

For a general product,
\begin{gather*}
{\rm Ric}_{\alpha \beta} = {\rm Ric}^{(1)}_{\alpha \beta},\qquad
{\rm Ric}_{jk} = {\rm Ric}^{(2)}_{jk},\qquad
{\rm Ric}_{\alpha k} = 0
\end{gather*}
and
\[
S = S^{(1)} + S^{(2)}.
\]
In particular, for $1 \times 3$ products ${\rm Ric}_{00} = 0$ and $S = S^{(2)}$.

The Bach tensor is somewhat more complicated. While the Bach tensor splits in the sense that the components with mixed indices are zero, the components corresponding to one factor depend on the curvature from the other factor.

The $1 \times 3$ and $2 \times 2$ cases are as follows: If $\dim\big(N^{(1)}\big) = 1$ and $\dim\big(N^{(2)}\big) = 3$ then
\begin{gather} \label{Bach1300eqn}
B_{00} = \left(-\frac{1}{12}\big(\Delta^{(2)} S^{(2)}\big)
		- \frac{1}{4}\left[ \big(|{\rm Ric}|^{(2)}\big)^2 - \frac{1}{3}\big(S^{(2)}\big)^2 \right] \right) g_{00},\\
B_{jk} = \frac{1}{2} \Delta^{(2)} {\rm Ric}^{(2)}_{jk}
				- \frac{1}{12} \Delta^{(2)} S^{(2)} g_{jk}
				- \frac{1}{6} S^{(2)}_{;jk}
			 - 2 \tr^{(2)}\big({\rm Ric}^{(2)} \otimes {\rm Ric}^{(2)}\big)_{jk}\nonumber\\
\hphantom{B_{jk} =}{}
				+ \frac{7}{6} S^{(2)} {\rm Ric}^{(2)}_{jk}
+ \frac{3}{4} \big(|{\rm Ric}|^{(2)}\big)^2 g_{jk}
				- \frac{5}{12} \big(S^{(2)}\big)^2 g_{jk},\label{Bach13123eqn}
\end{gather}
and
\[
B_{0 k} = 0.
\]
Here, $\tr({\rm Ric} \otimes {\rm Ric})_{jk} = g^{il} {\rm Ric}_{ij} {\rm Ric}_{lk}$. If $\dim\big(N^{(1)}\big) = \dim\big(N^{(2)}\big) = 2$ then
\begin{gather} \label{Bach22eqn1}
B_{\alpha \beta} = -\frac{1}{6} S^{(1)}_{;\alpha \beta}
		+ \frac{1}{6} \left[\Delta^{(1)} S^{(1)} - \frac{1}{2} \Delta^{(2)} S^{(2)}
						+ \frac{1}{4} \big(\big(S^{(1)}\big)^2 - \big(S^{(2)}\big)^2 \big) \right] g_{\alpha \beta}
\end{gather}
and similarly
\begin{gather} \label{Bach22eqn2}
B_{jk} = -\frac{1}{6} S^{(2)}_{;jk}
		+ \frac{1}{6} \left[\Delta^{(2)} S^{(2)} - \frac{1}{2} \Delta^{(1)} S^{(1)}
								+ \frac{1}{4} \big(\big(S^{(2)}\big)^2 - \big(S^{(1)}\big)^2 \big) \right] g_{jk}
\end{gather}
and
\[
B_{\alpha k} = 0.
\]

Using the formulation for the Bach tensor in the $1 \times 3$ setting, we have the following:

\begin{Proposition} \label{EinsteinSliceprop}
Let $\dim\big(N^{(1)}\big) = 1$ and $\dim\big(N^{(2)}\big) = 3$, and suppose $S = S^{(2)}$ is constant. Then $B = 0$ if and only if $g^{(2)}$ is Einstein.
\end{Proposition}

\begin{proof}If $g^{(2)}$ is Einstein then ${\rm Ric}^{(2)} = \frac{S^{(2)}}{3} g^{(2)}$ and as a result, equations~\eqref{Bach1300eqn} and~\eqref{Bach13123eqn} both simplify to zero.

In the other direction, if $g^{(2)}$ is not Einstein, then
\[
\big(|{\rm Ric}|^{(2)}\big)^2 > \frac{\big(S^{(2)}\big)^2}{3}.
\]
Since $\Delta^{(2)} S^{(2)} = 0$, then in particular, $B_{00} < 0$.
\end{proof}

An immediate consequence of this result and its proof is the following:

\begin{Proposition} \label{g00monotonicprop}
Let $M$ be a $1 \times 3$ product with product metric~$g$ solving equation~\eqref{bachfloweqn}. Suppose that for all time, the scalar curvature is constant on~$M$. Then $g_{00}$ is static if and only if $B = 0$, in which case, all components of the metric are static. Otherwise, $g_{00}$ is strictly decreasing.
\end{Proposition}

\section{Locally homogeneous spaces and Lie groups} \label{spacebackgroundsec}

A Riemannian manifold $(M, g)$ is locally homogeneous if for all points $p$, $q$ in $M$, there exist neighborhoods $U$ and $V$ about $p$ and $q$ respectively, and an isometry
\[
\varphi\colon \ U \rightarrow V
\]
with $\varphi(p) = q$. If, for all pairs of points, the isometry can be chosen to be global, so that $\varphi\colon M \rightarrow M$, then $(M, g)$ is homogeneous. If $M$ is closed and locally homogeneous, then its universal cover is homogeneous. A straightforward, but useful, result is that if a manifold is locally homogeneous, its scalar curvature is constant. There are nine three-dimensional simply connected homogeneous manifolds with compact quotients, six of which are Lie groups. The Lie groups support a larger class of homogeneous metrics and require a more sophisticated analysis than the three non-Lie groups. See \cite{IJ} for more details surrounding these definitions and results.

\subsection{Structure constants and curvature}

Let $G$ be a Lie group with Lie algebra $\mathfrak{g}$, and let $\{e_i\}$ be a left-invariant basis for $\mathfrak{g}$. The bracket can be expressed in terms of structure constants ${C_{ij}}^k$
\[
[e_i, e_j] = {C_{ij}}^k e_k.
\]

Given a left-invariant metric $g$, and working with a left-invariant frame, covariant derivatives, and then curvature can be expressed in terms of structure constants. The Ricci and scalar curvatures are
\begin{gather}\label{Riccistructureeqn}
{\rm Ric}_{jk} = -\frac{1}{2} \big({{C^l}_j}^p + {{C^p}_j}^l\big) C_{lkp}
			+ \frac{1}{4} {C^{lp}}_j C_{lpk}
			+ \frac{1}{2} {C^{lp}}_l(C_{pjk} + C_{pkj})
\end{gather}
and
\begin{gather} \label{Scalarstructureeqn}
S = -\frac{1}{4} C^{lkp} C_{lkp} - \frac{1}{2} C^{pkl} C_{lkp} - {C^{lp}}_l {C^k}_{pk}.
\end{gather}
Additionally, with a bit more calculation, the Laplacian of a left-invariant tensor can be expressed in terms of structure constants. For this paper, the Laplacian of a left-invariant symmetric $\binom{2}{0}$-tensor $T$ is needed and we have the following:
\begin{gather}
(\Delta T)_{ij} = \frac{1}{2} T_{pq}
				\big({{C^k}_i}^p {C_{kj}}^q + {C^{kp}}_i {{C_k}^q}_j + {C_i}^{pk} {{C_j}^q}_k
				-
					{{C^k}_i}^p {{C_k}^q}_j - {{C^k}_j}^p {{C_k}^q}_i \nonumber\\
\hphantom{(\Delta T)_{ij} = \frac{1}{2} T_{pq}\big(}{} -
					{{C^k}_i}^p {{C_j}^q}_k - {{C^k}_j}^p {{C_i}^q}_k
				 +
					{C^{kp}}_i {{C_j}^q}_k + {C^{kp}}_j {{C_i}^q}_k \big)\nonumber \\
\hphantom{(\Delta T)_{ij} =}{}+ \frac{1}{4} T_{qj}
				\big(\big({C^{kp}}_i - {{C^k}_i}^p + {C_i}^{pk}\big)
						\big({{C_k}^q}_p - {C_{kp}}^q + {{C_p}^q}_k\big)
 +2
					{C^{kp}}_k\big({{C_p}^q}_i - {C_{pi}}^q\big) \big) \nonumber\\
\hphantom{(\Delta T)_{ij} =}{}+ \frac{1}{4} T_{qi}
				\big(\big({C^{kp}}_j\! - {{C^k}_j}^p \!+ {C_j}^{pk}\big)
						\big({{C_k}^q}_p\! - {C_{kp}}^q\! + {{C_p}^q}_k\big)
 +2
					{C^{kp}}_k\big({{C_p}^q}_j \! - {C_{pj}}^q\big) \big).\!\!\!\!\!\!\label{Laplacianstructureeqn}
\end{gather}

\subsection{Three-dimensional Lie groups}

As seen in \cite{RyanShepley}, the six three-dimensional simply connected Lie groups with compact quotients are all unimodular and all have the property that for each group, there is a basis for the Lie algebra such that the structure constants can be represented by
\[
{C_{ij}}^k = \varepsilon_{ijs} E^{ks},
\]
where $\varepsilon_{ijk}$ is the Levi-Civita symbol which captures the parity of the permutation generating ``$ijk$'' with $\varepsilon_{123} = 1$, and where~$E$ is a $3 \times 3$ matrix specific to the group. See Fig.~\ref{bianchifigure}.

\begin{figure}\centering
\begin{tabular}{c|c|l}
Bianchi & \multirow{2}{*}{Group} & \multirow{2}{*}{$\ \ \ \ \ \ E$} \\
Type & & \\
\hline
\hline
I & $\mathbb{R}^3$ & 0 \tsep{2pt}\\
\hline
II & ${\rm Nil}$ & $\diag(1, 0, 0)$\tsep{1pt} \\
\hline
VI$_{0}$ & ${\rm Solv}$ & $\diag(-1, 1, 0)$\tsep{1pt} \\
\hline
VII$_0$ & ${\rm E}(2)$ & $\diag(-1, -1, 0)$\tsep{1pt} \\
 \hline
VIII & $\widehat{{\rm SL}}(2, \mathbb{R})$ & $\diag(-1, 1, 1)$ \tsep{4pt}\\
\hline
IX & $\mathbb{S}^3$ & ${\rm id}$ \tsep{2pt}
\end{tabular}
\caption{The six three-dimensional simply connected Lie groups with compact quotients.
The $3 \times 3$ matrix $E$ encapsulates the structure constants.} \label{bianchifigure}
\end{figure}

To simplify the later analysis, the Bach flow equations will be determined in a basis where the structure constants have the form indicated here and where the initial metric is diagonal. As shown in \cite{Milnor}, such an initial set-up is always possible:

\begin{Theorem} \label{metricstructureconstantsthm}
Given a three-dimensional Lie algebra with structure constants of the form
\begin{gather} \label{structureconstmeqn}
{C_{ij}}^k = \varepsilon_{ijl} E^{lk},
\end{gather}
and an inner product $g$, there is a basis where
\begin{itemize}\itemsep=0pt
\item $g$ is diagonal $($the basis is orthogonal$)$,
\item the structure constants can still be written in the form \eqref{structureconstmeqn},
\item the matrix $E$ is unchanged.
\end{itemize}
\end{Theorem}

The proof of this theorem follows from the principal axis theorem and the fact that structure constants can be rescaled by rescaling the basis. We call the basis guaranteed by Theorem~\ref{metricstructureconstantsthm} a diagonalizing basis.

In light of the structure afforded by Theorem \ref{metricstructureconstantsthm}, we note the following general facts about the Ricci tensor and Bach tensor:

\begin{Proposition} \label{RicBachdiagprop}
Let $N$ be a three-dimensional Lie group with structure constants of the form of equation \eqref{structureconstmeqn}, and left-invariant metric $g$. Then, in a basis where $E$ and $g$ are diagonal, the Ricci tensor is diagonal and, on $\mathbb{S}^1 \times N$, the Bach tensor is diagonal.
\end{Proposition}

For the Ricci tensor, this was established in \cite{IJ} for the specific matrices in Fig.~\ref{bianchifigure}. This proof shows that the property is a consequence of the diagonal structure of~$E$ and~$g$, and not specific to particular matrices. The proof of the general result follows from careful accounting of the indices in each term found in the formulas for the Ricci and Bach tensors, using the fact that in three dimensions, the indices are restricted to just three values.

\begin{proof}First observe that if $E$ is diagonal then equation \eqref{structureconstmeqn} shows that ${C_{ij}}^k$ can only be nonzero if $i$, $j$, and $k$ are all different. Moreover, if $g$ is diagonal, then the same must be true for any raising or lowering of any of the indices. Hence any structure constant with a repeated index must be zero, and in any double sum involving a pair of structure constants, the two free indices must be equal in order for the result to be nonzero. Based on these observations, every term in equation~\eqref{Riccistructureeqn} must be zero unless $j = k$ so ${\rm Ric}$ must be diagonal.

The analysis for the Bach tensor is similar. Looking at equation~\eqref{Bach13123eqn}, note first that the second and third terms are zero since scalar curvature is constant, and the fifth, sixth, and seventh terms are diagonal since ${\rm Ric}^{(2)}$ and $g$ are diagonal. So the only terms to check are the first and the fourth. For the fourth term, we have
\[
\tr({\rm Ric} \otimes {\rm Ric})_{jk} = g^{il} {\rm Ric}_{ij} {\rm Ric}_{lk}.
\]
Since $g$ is diagonal, the terms in this sum are only nonzero when $i = l$, and then, since ${\rm Ric}$ is diagonal, we can only have a nonzero term when $j = k$.

Finally for the first term, we use equation \eqref{Laplacianstructureeqn} with $T = {\rm Ric}^{(2)}$ to analyze $\big(\Delta^{(2)} {\rm Ric}^{(2)}\big)_{ij}$. Equation \eqref{Laplacianstructureeqn} has three large terms in it. For the first term, since $T$ is diagonal, the only way any of the sums of products can be nonzero is if $p = q$, but then each product becomes a double sum and so must be zero unless $i = j$. For the second and third terms, one piece is zero because of a structure constant with a repeated index. For the rest, the double sums again require the third pair of indices to match in order to produce something nonzero, and since~$T$ is diagonal, the only nonzero terms appear when $i = j$.
\end{proof}

To help with the analysis of curvature along the flow, we have the following:

\begin{Lemma} \label{sectionalcurvaturelemma}Let $\{e_1, e_2, e_3\}$ be an orthogonal basis for the tangent space of a point in a $3$-dimensional manifold. Then at that point, the sectional curvatures are given by
\[
K(e_i, e_j) = \frac{{\rm Ric}_{ii}}{g_{ii}} + \frac{{\rm Ric}_{jj}}{g_{jj}} - \frac{S}{2}.
\]
\end{Lemma}

\begin{proof}
On a 3-dimensional manifold, Riemann curvature can be expressed completely in terms of Ricci and scalar curvature as
\[
R = {\rm Ric} \circ g - \frac{S}{4} g \circ g,
\]
where $A \circ B$ is the Kulkarni--Nomizu product. Using an orthogonal basis, this reduces to
\[
R_{ijji} = {\rm Ric}_{ii} g_{jj} + {\rm Ric}_{jj} g_{ii} - \frac{S}{2} g_{ii} g_{jj}.
\]
Sectional curvature is given by
\[
K(v, w) = \frac{R(v, w, w, v)}{|v|^2 |w|^2 - \langle v, w \rangle^2},
\]
so
\begin{gather*}
K(e_i, e_j) = \frac{R(e_i, e_j, e_j, e_i)}{|e_i|^2 |e_j|^2 - \langle e_i, e_j \rangle^2}
		 = \frac{{\rm Ric}_{ii} g_{jj} + {\rm Ric}_{jj} g_{ii} - \frac{S}{2} g_{ii} g_{jj}}{g_{ii} g_{jj}} = \frac{{\rm Ric}_{ii}}{g_{ii}} + \frac{{\rm Ric}_{jj}}{g_{jj}} - \frac{S}{2}
\end{gather*}
as desired.
\end{proof}

Lemma \ref{sectionalcurvaturelemma} is all that is needed in this paper since the four-dimensional manifolds considered are $1 \times 3$ products, so the formula above can be used for the three-dimensional slice, and the sectional curvatures involving the one-dimensional slice are zero.

\section{Ordinary differential equations} \label{odesec}

The ordinary differential equations to which Bach flow reduces on homogeneous products are analyzed using standard techniques which are recalled here. First, we appeal to existence and uniqueness of solutions regularly and without mention. In some instances, the equations of interest are separable and explicit solutions may be found. When such explicit solutions cannot be found, the Escape Lemma, which states that if a maximal flow does not exist for all time then it cannot lie in a compact set, may be used to help determine the qualitative behavior of solutions. See~\cite{Jack2} for details surrounding these results. In addition to these methods, we make use of a couple more specialized results which follow.

The following technical lemma provides a technique for determining the images of the integral curves for a system of two equations involving homogeneous functions.

\begin{Lemma} \label{homogeneousodesolnlemma}
Let $(x(t), y(t))$ solve the following system
\begin{equation*}
\frac{{\rm d}x}{{\rm d}t} = p(x, y),\qquad \frac{{\rm d}y}{{\rm d}t} = q(x, y),
\end{equation*}
where $p$ and $q$ are both homogeneous of degree $k$. Suppose $x \neq 0$, $p(x,y) \neq 0$, and $\frac{q(x,y)}{p(x,y)} \neq \frac{y}{x}$. Then $(x(t), y(t))$ will lie in the curve
\[
x = \eta e^{\Psi \left(\frac{y}{x}\right)},
\]
where
\[
\Psi(v) = \int \frac{1}{\frac{q(1, v)}{p(1, v)} - v} {\rm d}v
\]
and $\eta$ is a constant depending on the initial conditions.
\end{Lemma}

To prove this result, briefly, express $\frac{{\rm d}y}{{\rm d}x}$ in terms of $p$ and $q$ and then compute the derivative of $v = \frac{y}{x}$ with respect to~$x$, substitute and use separation of variables. The details are left to the reader.

In general, an integral curve for a vector field can be bounded but fail to converge to a limit. The next lemma shows that if a coordinate of a bounded solution does converge, then that component of the vector field must go to zero. The proof is left to the reader.

\begin{Lemma} \label{odeconvergencelemma}
Let $x(t) = \big(x^1(t), \ldots, x^n(t)\big)$ be a bounded solution to
\[
\dtime x = V(x),
\]
where $V$ is a continuous vector field on a domain $D$. Suppose
\[
\lim_{t \rightarrow \infty} x^i(t) = L
\]
and let $\{x_k\} = \{x(t_k)\}$, $t_k \rightarrow \infty$ be a sequence of points on the curve that converges to $x_{\infty} \in D$. Then $V^i(x_{\infty}) = 0$.
\end{Lemma}

\section[Bach flow on locally homogeneous $1 \times 3$ products]{Bach flow on locally homogeneous $\boldsymbol{1 \times 3}$ products} \label{13flowsec}

In this section, the main results of this paper are proved for $1 \times 3$ products that are not also $2 \times 2$ products. For each universal cover, explicit formulas for the Bach tensor are found and the evolution of the metric under Bach flow is determined. In some cases, explicit solutions are found. When explicit solutions are not found, qualitative behavior is determined. Limiting behaviors of both the metric and its curvature are also found, and convergence, in the Gromov--Hausdorff or pointed Gromov--Hausdorff topology, is described.

The general method is as follows: Given an initial metric, a diagonalizing basis is found so that the metric is diagonal
\[
h = \diag(h_{00}, h_{11}, h_{22}, h_{33}).
\]
Its Ricci and scalar curvatures are calculated using equations \eqref{Riccistructureeqn} and \eqref{Scalarstructureeqn}, and then using equations \eqref{Bach1300eqn}, \eqref{Bach13123eqn}, and \eqref{Laplacianstructureeqn} the Bach tensors are calculated. As indicated by Proposition~\ref{RicBachdiagprop}, the Bach tensor is also diagonal, so the fact that the metric is diagonal is preserved along the flow. The solution will be denoted
\[
g = \diag(g_{00}, g_{11}, g_{22}, g_{33}).
\]

One quantity that makes a regular appearance is
\[
\beta = \frac{1}{6 (\det g)^2} = \frac{1}{6 (\det h)^2}.
\]
This quantity depends on the initial metric, but is constant along the flow. As a consequence, once an initial metric is chosen, $\beta$ can be treated as a constant for the whole system.

After explicitly determining the Bach tensor in each case, a general structure emerges. Speci\-fi\-cally, the nonzero components of the Bach tensor have the form
\[
B_{ii} = \alpha_i \beta p_{i}(g_{11}, g_{22}, g_{33}) (g_{00})^2 g_{ii},
\]
where $\alpha_i$ is a constant and $p_{i}$ is a homogeneous fourth degree polynomial.
This structure makes a great deal of qualitative analysis possible when explicit solutions are not found. The details of the analysis vary from space to space, although there are similarities when the spaces themselves have similar structure constant matrices $E$.

%
%
One general fact is that $g_{00}$ is decreasing, as indicated by Proposition \ref{g00monotonicprop}. This fact will not be explicitly included in the specific theorems for each space. Another general fact is that the flow is defined (at least) on the interval $[0, \infty)$. In cases where the flow remains bounded, this follows from the escape lemma. In cases where the flow does not remain bounded, this is discovered after the analysis of each flow is completed and follows from the work in~\cite{Lopez}, which shows that the maximal time is finite only if there is curvature blow-up, and the fact that in all of our cases, curvature remains bounded.

\subsection[$\mathbb{R}^3$]{$\boldsymbol{\mathbb{R}^3}$} \label{r3subsection}

For this manifold the matrix $E$ used to determine the structure constants in Theorem~\ref{metricstructureconstantsthm} is the zero matrix, so regardless of the initial metric, the structure constants are all zero. Hence the Ricci tensor, scalar curvature are zero, and on $\mathbb{R} \times \mathbb{R}^3$ the Bach tensor is zero and so the metric is static under Bach flow.

\subsection[${\rm Nil}$]{$\boldsymbol{{\rm Nil}}$}

For this manifold the matrix used to determine the structure constants in Theorem~\ref{metricstructureconstantsthm} is
\[
E = \diag(1, 0, 0).
\]
For any metric $g$, using a diagonalizing basis, the Ricci tensor is diagonal with
\begin{gather*}
{\rm Ric}_{11} = \frac{(g_{11})^2}{2 g_{22} g_{33}},\qquad
{\rm Ric}_{22} = -\frac{g_{11}}{2 g_{33}},\qquad
{\rm Ric}_{33} = -\frac{g_{11}}{2 g_{22}}
\end{gather*}
and scalar curvature is
\[
S = -\frac{g_{00} (g_{11})^2}{2 \det g} .
\]
The Bach tensor on $\mathbb{R} \times {\rm Nil}$ is diagonal with
\begin{alignat*}{3}
& B_{00} = -\beta (g_{00})^3 (g_{11})^4, \qquad &&
B_{11} = -5 \beta (g_{00})^2 (g_{11})^5,&\\
& B_{22} = 3 \beta (g_{00})^2 (g_{11})^4 g_{22}, \qquad &&
B_{33} = 3 \beta (g_{00})^2 (g_{11})^4 g_{33}.&
\end{alignat*}

With the Bach tensor in hand, we have the following theorem:

\begin{Theorem} On $\widehat{M} = \mathbb{R} \times {\rm Nil}$ the solutions to equation \eqref{bachfloweqn} in a diagonalizing basis for $h$ are given by
\begin{alignat*}{3}
& g_{00}(t) = \big(\gamma t + (h_{00})^{-22} \big)^{-\frac{1}{22}}, \qquad &&
g_{11}(t) = \alpha \big(\gamma t + (h_{00})^{-22} \big)^{-\frac{5}{22}},& \\
& g_{22}(t) = h_{22} (h_{00})^{3}
		\big(\gamma t + (h_{00})^{-22} \big)^{\frac{3}{22}}, \qquad &&
g_{33}(t) = h_{33} (h_{00})^{3}
		\big(\gamma t + (h_{00})^{-22} \big)^{\frac{3}{22}},&
\end{alignat*}
where
\[
\alpha = \frac{h_{11}}{(h_{00})^5} \qquad \text{and} \qquad
\gamma = 22 \alpha^4 \beta
		= \frac{11}{3 (\det h)^2} \left(\frac{h_{11}}{(h_{00})^5}\right)^4.
\]
\end{Theorem}

\begin{proof}
Note that the first and second equations are coupled and the third and fourth equations depend on the first and second solutions, but are otherwise uncoupled. Because everything is multiplicative, we can solve explicitly for $g_{00}$, and $g_{11}$, and then $g_{22}$ and $g_{33}$.

Starting with the following ansantz
\begin{equation*} 
g_{11} = \alpha (g_{00})^k
\end{equation*}
and then comparing the resulting differential equations produces
\[
g_{00}(t) = \big(\gamma t + (h_{00})^{-22}\big)^{-\frac{1}{22}}
\qquad \text{and}\qquad g_{11}(t) = \alpha \big(\gamma t + (h_{00})^{-22} \big)^{-\frac{5}{22}},
\]
where
\[
\gamma = 22 \alpha^4 \beta.
\]

Then we can solve for $g_{22}$ and $g_{33}$:
\begin{gather*}
g_{22}(t) = h_{22} (h_{00})^{3}
		\big(\gamma t + (h_{00})^{-22} \big)^{\frac{3}{22}}
\qquad \text{and}\qquad
g_{33}(t) = h_{33} (h_{00})^{3}
		\big(\gamma t + (h_{00})^{-22}\big)^{\frac{3}{22}}.\tag*{\qed}
\end{gather*}\renewcommand{\qed}{}
\end{proof}

With these solutions in hand, we find two dimensions collapse in the limit as $t \rightarrow \infty$. The~``$g_{00}$'' direction collapses more slowly than the first dimension in ${\rm Nil}$. Meanwhile, the other two dimensions grow at the same rate, preserving their aspect ratio. These solutions are immortal, but not ancient.

All components of the Ricci tensor converge to zero in the limit, and using Lemma~\ref{sectionalcurvaturelemma} we have the following:

\begin{Theorem}Let $M$ be a compact quotient of $\mathbb{R} \times {\rm Nil}$ and let $p \in M$. Let $g$ solve equation~\eqref{bachfloweqn} where $h$ is locally homogeneous. Then $(M, g, p)$ collapses to a flat surface in the pointed Gromov--Hausdorff topology.
\end{Theorem}

\subsection[${\rm Solv}$]{$\boldsymbol{{\rm Solv}}$} \label{solvsection}

For this manifold the matrix used to determine the structure constants in Theorem~\ref{metricstructureconstantsthm} is
\[
E = \diag(-1, 1, 0).
\]
For any metric $g$, using a diagonalizing basis, the Ricci tensor is diagonal with
\begin{gather*}
{\rm Ric}_{11} = \frac{(g_{11})^2 - (g_{22})^2}{2 g_{22} g_{33}},\qquad
{\rm Ric}_{22} = \frac{(g_{22})^2 - (g_{11})^2}{2 g_{11} g_{33}},\qquad
{\rm Ric}_{33} = - \frac{(g_{11} + g_{22})^2}{2 g_{11} g_{22}}
\end{gather*}
and scalar curvature is
\[
S = -\frac{(g_{11} + g_{22})^2}{2 g_{11} g_{22} g_{33}}.
\]
The Bach tensor is diagonal with
\begin{alignat*}{3}
& B_{00} = -\beta p_{\rm I}(g_{11}, g_{22}) (g_{00})^3, \qquad &&
B_{11} = -\beta p_{\rm II}(g_{11}, g_{22}) (g_{00})^2 g_{11}, & \\
& B_{22} = -\beta p_{\rm II}(g_{22}, g_{11}) (g_{00})^2 g_{22}, \qquad&&
B_{33} = 3\beta p_{\rm I}(g_{11}, g_{22}) (g_{00})^2 g_{33},&
\end{alignat*}
where
\begin{gather*}
p_{\rm I}(x, y) = x^4 + x^3y + xy^3 + y^4, \qquad
p_{\rm II}(x, y) = 5x^4 + 3x^3y - xy^3 - 3y^4.
\end{gather*}

With the Bach tensor in hand, we have the following theorem:

\begin{Theorem} \label{solvthm} On $\widehat{M} = \mathbb{R} \times {\rm Solv}$ every solution to equation~\eqref{bachfloweqn} in a diagonalizing basis has the following properties:
\begin{itemize}\itemsep=0pt
\item $g_{00}, g_{11}, g_{22} \rightarrow 0$;
\item $g_{33} \rightarrow \infty$ monotonically;
\item $\frac{g_{11}}{g_{22}} \rightarrow 1$.
\end{itemize}
If $h_{11} = h_{22}$, then
\begin{gather*}
g_{00} = \mu^{\frac{1}{2}} \bigl(24 \mu \beta t + (h_{11})^{-6} \bigr)^{-\frac{1}{6}},\qquad
g_{11} = g_{22} = \bigl(24 \mu \beta t + (h_{11})^{-6} \bigr)^{-\frac{1}{6}}, \\
g_{33} = (h_{11})^3 h_{33} \bigl(24 \mu \beta t + (h_{11})^{-6} \bigr)^{\frac{1}{2}}.
\end{gather*}
Otherwise, if (without loss of generality) $h_{11} < h_{22}$, then
\begin{itemize}\itemsep=0pt
\item $g_{11} < g_{22}$ for the entire flow;
\item $g_{22}$ is decreasing;
\item $\frac{g_{11}}{g_{22}}$ is increasing;
\item $g_{11}$ and $g_{22}$ are related by
\[
(g_{11} g_{22})^{25} = \eta (g_{22}-g_{11})^{4} \big(2(g_{22})^2 + g_{11} g_{22} + 2(g_{11})^2\big)^{3},
\]
where
\[
\eta = \frac{(h_{11} h_{22})^{25}}{(h_{22}-h_{11})^{4}
				\big(2(h_{22})^2 + h_{11} h_{22} + 2(h_{11})^2\big)^{3}}.
\]
\end{itemize}
\end{Theorem}

It turns out that in addition to appearing in the Bach tensor above, the two polynomials~$p_{\rm I}$ and~$p_{\rm II}$ also make an appearance in the next section so we establish some facts about them for use here and later.

\begin{Lemma}The polynomial $p_{\rm I}(x,y)$ is symmetric, homogeneous of degree~$4$, positive when~$x$ or~$y$ is nonzero, and can be factored as
\[
p_{\rm I}(x, y) = (x+y)^2\big(x^2 - xy + y^2\big).
\]
The polynomial $p_{\rm II}(x,y)$ is homogeneous of degree $4$ and can be factored as
\[
p_{\rm II}(x, y) = (x+y)\big(5x^3 - 2x^2y + 2xy^2 - 3y^3\big).
\]
The cubic factor has exactly one real factor $(\alpha x - y)$ where $\alpha$ is about~$1.23$. If $x > y$, $p_{\rm II}(x, y) > 0$.
\end{Lemma}

The proof of this lemma is left to the reader. While $\alpha$, as the root of a cubic, can be found exactly, this exact form is not important for the analysis here. With these facts about~$p_{\rm I}$ and~$p_{\rm II}$ established, we proceed with the proof of Theorem~\ref{solvthm}.

\begin{proof}[Proof of Theorem \ref{solvthm}]
Since $B_{00}$ and $B_{33}$ are so similar, we can compute
\begin{equation*}
\frac{\dtime g_{33}}{\dtime g_{00}} = \frac{B_{33}}{B_{00}}
			= \frac{-3 g_{33}}{g_{00}},
\end{equation*}
which implies
\begin{gather} \label{solvg33eqn}
g_{33} = \gamma (g_{00})^{-3},
\end{gather}
where $\gamma = (h_{00})^3 h_{33}$. Since, by Proposition~\ref{g00monotonicprop}, $g_{00}$ is decreasing, this shows $g_{33}$ must be increasing. Since $\det g$ is constant, using equation \eqref{solvg33eqn} we have
\begin{gather} \label{solvg00eqn}
(g_{00})^2 = \mu g_{11} g_{22},
\end{gather}
where
\[
\mu = \frac{\gamma}{\det h}.
\]
Incorporating these identities into the formulas for $B_{11}$ and $B_{22}$ we have
\begin{gather} \label{solvrestrictedg11eqn}
\dtime g_{11} = -\mu \beta p_{\rm II}(g_{11}, g_{22}) (g_{11})^2 g_{22}
\end{gather}
and
\[
\dtime g_{22} = -\mu \beta p_{\rm II}(g_{22}, g_{11}) (g_{22})^2 g_{11}.
\]

Because of the symmetry in these equations, we may, without loss of generality, restrict our attention to the region defined by $0 \leq g_{11} \leq g_{22}$.

If $h_{11} = h_{22}$ then, focusing on $g_{11}$ we get
\[
\dtime g_{11} = -4\mu \beta g_{11}^7,
\]
which is separable. Solving, we get
\[
g_{11} = g_{22} = \big(24 \mu \beta t + (h_{11})^{-6} \big)^{-\frac{1}{6}}.
\]

Next, we solve for $g_{00}$ and $g_{33}$ using equations \eqref{solvg33eqn} and \eqref{solvg00eqn}. We have
\[
g_{00} = \mu^{\frac{1}{2}} \bigl(24 \mu \beta t + (h_{11})^{-6} \bigr)^{-\frac{1}{6}}
\qquad \text{and} \qquad
g_{33} = (h_{11})^3 h_{33} \bigl(24 \mu \beta t + (h_{11})^{-6} \bigr)^{\frac{1}{2}}.
\]
In this special case, we see that under Bach flow, any compact quotient of $\mathbb{R} \times {\rm Solv}$ collapses to a curve in the limit.

If $h_{11} < h_{22}$ then $g_{11} < g_{22}$ by existence and uniqueness, since we have a solution that preserves the equality $g_{11} = g_{22}$. With this inequality preserved, from the properties of $p_{\rm II}$, looking at equation~\eqref{solvrestrictedg11eqn} we find that $g_{22}$ is decreasing. Also
\begin{align*}
\dtime \left(\frac{g_{11}}{g_{22}} \right)
	&= \frac{\big(\dtime g_{11} \big) g_{22} - g_{11} \big(\dtime g_{22} \big)}{(g_{22})^2} \\
	&= 4 \mu \beta (g_{11} + g_{22})(g_{22} - g_{11})
					\big(2(g_{11})^2 + g_{11} g_{22} + 2 (g_{22})^2 \big) (g_{11})^2 > 0
\end{align*}
so we find that $\frac{g_{11}}{g_{22}}$ is increasing. This fraction is bounded above by 1, so it must converge. We will see below that it converges to~1.

The fact that this fraction is increasing also implies that $g_{11}$ cannot converge to zero unless~$g_{22}$ does as well. This, combined with Lemma~\ref{odeconvergencelemma} implies that $g_{22}$ and hence $g_{11}$ must converge to~0 since the only points in the domain of interest where $\dtime g_{22}$ is zero are along the~$g_{22}$ axis.

Next, letting $v = \frac{g_{11}}{g_{22}}$ we have
\begin{gather*}
\frac{{\rm d} g_{11}}{{\rm d} g_{22}} = \frac{p_{\rm II}(v, 1) v^2}{p_{\rm II}(1, v) v}
		 = v \frac{\big(5v^3 - 2v^2 + 2v- 3\big)}{\big(5 - 2v + 2v^2 - 3v^3\big)}.
\end{gather*}
Therefore by Lemma~\ref{homogeneousodesolnlemma}, the solution curves for our original differential equation satisfy the equation
\[
g_{22} = \frac{\tilde{\eta} \big(1-\frac{g_{11}}{g_{22}}\big)^{\frac{1}{10}}
		\left(2 + \frac{g_{11}}{g_{22}} + 2\big(\frac{g_{11}}{g_{22}}\big)^2\right)^{\frac{3}{40}}}				{\big(\frac{g_{11}}{g_{22}}\big)^{\frac{5}{8}}},
\]
where $\tilde{\eta}$ is a constant determined by the initial conditions.
Multiplying both sides by $\big(\frac{g_{11}}{g_{22}}\big)^{\frac{5}{8}}$, we have
\begin{gather} \label{solvcurveeqn}
\left(\frac{g_{11}}{g_{22}}\right)^{\frac{5}{8}} g_{22}
		= \tilde{\eta} \left(1-\frac{g_{11}}{g_{22}}\right)^{\frac{1}{10}}
		\left(2 + \frac{g_{11}}{g_{22}} + 2\left(\frac{g_{11}}{g_{22}}\right)^2\right)^{\frac{3}{40}}.
\end{gather}
This can be rewritten as
\begin{equation*}
(g_{11} g_{22})^{25} = \eta (g_{22}-g_{11})^{4} \big(2(g_{22})^2 + g_{11} g_{22} + 2(g_{11})^2\big)^{3},
\end{equation*}
where $\eta = \tilde{\eta}^{40}$. This is true in particular at $t = 0$, so
\[
\eta = \frac{(h_{11} h_{22})^{25}}{(h_{22}-h_{11})^{4}
				\big(2(h_{22})^2 + h_{11} h_{22} + 2(h_{11})^2\big)^{3}}.
\]
Taking the limit as $t \rightarrow \infty$, the left side of equation \eqref{solvcurveeqn} must be zero, and therefore so must the right. Since the second factor is positive, it follows that $\lim\limits_{t \rightarrow \infty}\big(1-\frac{g_{11}}{g_{22}}\big)^{\frac{1}{10}} = 0$ and so $\lim\limits_{t \rightarrow \infty} \frac{g_{11}}{g_{22}} = 1$.

From the analysis above, we know that in general, $g_{11}$ and $g_{22}$ go to zero as $t$ goes to infinity. From equations \eqref{solvg00eqn} and \eqref{solvg33eqn}, we then know that $g_{00}$ also goes to zero and $g_{33}$ grows to infinity. Therefore, under Bach flow, any compact quotient of $\mathbb{R} \times {\rm Solv}$ collapses to a curve in the limit. These facts, combined with the fact that $\frac{g_{11}}{g_{22}}$ goes to 1 imply that a general solution approaches the specific solution found above in the limit.
\end{proof}

We now have the following:

\begin{Theorem}Let $M$ be a compact quotient of $\mathbb{R} \times {\rm Solv}$ and let $p \in M$. Let $g$ solve equa\-tion~\eqref{bachfloweqn} where $h$ is locally homogeneous. Then $(M, g, p)$ collapses to a line in the pointed Gromov--Hausdorff topology.
\end{Theorem}

\begin{proof}From the previous theorem, we know that three dimensions collapse, while one expands. Moreover, working in a diagonalizing basis, we find that ${\rm Ric}_{11}$ and ${\rm Ric}_{22}$ converge to zero while ${\rm Ric}_{33}$ converges to $-2$. Therefore, by Lemma~\ref{sectionalcurvaturelemma}, all the sectional curvatures go to zero in the limit.
\end{proof}

\subsection[${\rm E}(2)$]{$\boldsymbol{{\rm E}(2)}$}

For this manifold the matrix used to determine the structure constants in Theorem~\ref{metricstructureconstantsthm} is
\[
E = \diag(-1, -1, 0).
\]
For any metric $g$, using a diagonalizing basis, the Ricci tensor is diagonal with
\begin{gather*}
{\rm Ric}_{11} = \frac{(g_{11})^2 - (g_{22})^2}{2 g_{22} g_{33}},\qquad
{\rm Ric}_{22} = \frac{(g_{22})^2 - (g_{11})^2}{2 g_{11} g_{33}},\qquad
{\rm Ric}_{33} = - \frac{(g_{11} - g_{22})^2}{2 g_{11} g_{22}}
\end{gather*}
and scalar curvature is
\[
S = -\frac{(g_{11} - g_{22})^2}{2 g_{11} g_{22} g_{33}}.
\]
The Bach tensor is diagonal with
\begin{alignat*}{3}
& B_{00}= -\beta p_{\rm I}(-g_{11}, g_{22}) (g_{00})^3,\qquad &&
B_{11}= -\beta p_{\rm II}(-g_{11}, g_{22}) (g_{00})^2 g_{11}, &\\
& B_{22}= -\beta p_{\rm II}(g_{22}, -g_{11}) (g_{00})^2 g_{22}, \qquad &&
B_{33}= 3\beta p_{\rm I}(-g_{11}, g_{22}) (g_{00})^2 g_{33},&
\end{alignat*}
where $p_{\rm I}$ and $p_{\rm II}$ were defined is Section \ref{solvsection}.

\begin{Theorem} \label{e2thm}
On $\widehat{M} = \mathbb{R} \times {\rm E}(2)$ every solution to equation \eqref{bachfloweqn} in a diagonalizing basis has the following properties:
\begin{itemize}\itemsep=0pt
\item $g_{11}$ and $g_{22}$ are related by
\[
(g_{11} g_{22})^{25} = \eta (g_{22} + g_{11})^{4} \big(2(g_{22})^2 - g_{11} g_{22} + 2(g_{11})^2\big)^{3},
\]
where
\[
\eta = \frac{(h_{11} h_{22})^{25}}{(h_{22} + h_{11})^{4}
				\big(2(h_{22})^2 - h_{11} h_{22} + 2(h_{11})^2\big)^{3}};
\]
\item the flow exists for all time and as $t \rightarrow \infty$,
\begin{gather*}
g_{11}, g_{22} \rightarrow (432 \eta)^{\frac{1}{40}}, \qquad
g_{00} \rightarrow \left(\frac{(h_{00})^3 h_{33}}{\det h} \right)^{\frac{1}{2}} (432 \eta)^{\frac{1}{40}}, \\
g_{33} \rightarrow \left( \frac{(\det h)^{3}}{(h_{00})^{3} h_{33}} \right)^{\frac{1}{2}} (432 \eta)^{-\frac{3}{40}};
\end{gather*}
\item $g_{33}$ is increasing.
\end{itemize}
If $h_{11} = h_{22}$, then the solution is static. Otherwise, if (without loss of generality) $h_{11} < h_{22}$, then
\begin{itemize}\itemsep=0pt
\item $g_{11} < g_{22}$ for the entire flow,
\item $g_{11}$ is increasing,
\item $g_{22}$ is decreasing.
\end{itemize}
\end{Theorem}

\begin{proof}
Since the only difference between this system and that of $\mathbb{R} \times {\rm Solv}$ is the minus sign on one of the variables in $p_{\rm I}$ and $p_{\rm II}$, most of the initial analysis of the previous section carries over and we have
\begin{gather} \label{e2g33eqn}
g_{33} = \gamma (g_{00})^{-3}.
\end{gather}
with $\gamma = (h_{00})^3 h_{33}$. Because of this inverse relationship, since $g_{00}$ is decreasing, we find that $g_{33}$ must be increasing. We also have
\begin{gather} \label{e2g00eqn}
(g_{00})^2 = \mu g_{11} g_{22},
\end{gather}
where
\[
\mu = \frac{\gamma}{\det h}.
\]

Incorporating these identities into the formulas for $B_{11}$ and $B_{22}$ we have
\begin{gather*}
\dtime g_{11} = -\mu \beta p_{\rm II}(-g_{11}, g_{22}) (g_{11})^2 g_{22}, \qquad
\dtime g_{22} = -\mu \beta p_{\rm II}(g_{22}, -g_{11}) (g_{22})^2 g_{11}.
\end{gather*}

Because of the symmetry in these equations, we may, without loss of generality, restrict our attention to the region defined by $0 \leq g_{11} \leq g_{22}$. If $h_{11} = h_{22}$, then $B = 0$ and we have a set of stationary solutions corresponding to the flat metrics on $E(2)$. If $h_{11} < h_{22}$ then $g_{11} < g_{22}$ for all time and from the poperties of $p_{\rm II}$, we find that $g_{11}$ is increasing and $g_{22}$ is decreasing. Therefore, both must converge and by Lemma \ref{odeconvergencelemma} this can only happen at a point where $g_{11} = g_{22}$.

As with ${\rm Solv}$, we can say a bit more about the curves traced out by the solutions using Lemma \ref{homogeneousodesolnlemma}. Except for two minus signs, the analysis here is almost identical to that for ${\rm Solv}$ and we find
\begin{gather} \label{e2curveeqn}
(g_{11} g_{22})^{25} = \eta (g_{22}+g_{11})^{4} \big(2(g_{22})^2 - g_{11} g_{22} + 2(g_{11})^2\big)^{3},
\end{gather}
where
\[
\eta = \frac{(h_{11} h_{22})^{25}}{(h_{22}+h_{11})^{4}
				\big(2(h_{22})^2 - h_{11} h_{22} + 2(h_{11})^2\big)^{3}}.
\]

Let $g_{ii}(\infty)$ be the limit of $g_{ii}$ as $t \rightarrow \infty$. Then we know that $g_{11}(\infty) = g_{22}(\infty)$ and using equation \eqref{e2curveeqn} we find
\[
\bigl(g_{11}(\infty)\bigr)^{50} = \eta \bigl(2g_{11}(\infty)\bigr)^4 \bigl(3 (g_{11}(\infty))^2 \bigr)^3
\]
and so
\[
g_{11}(\infty) = g_{22}(\infty) = (432 \eta)^{\frac{1}{40}}.
\]
Then by equation \eqref{e2g00eqn}
\[
g_{00}(\infty) = \mu^{\frac{1}{2}} (432 \eta)^{\frac{1}{40}}
			= \left(\frac{(h_{00})^3 h_{33}}{\det h} \right)^{\frac{1}{2}} (432 \eta)^{\frac{1}{40}}
\]
and by equation \eqref{e2g33eqn}
\begin{gather*}
g_{33}(\infty) = \gamma \bigl[\mu^{\frac{1}{2}} (432 \eta)^{\frac{1}{40}}\bigr]^{-3}
			= \left( \frac{(\det h)^{3}}{(h_{00})^{3} h_{33}} \right)^{\frac{1}{2}} (432 \eta)^{-\frac{3}{40}}.\tag*{\qed}
\end{gather*} \renewcommand{\qed}{}
\end{proof}

Finally, we have the following

\begin{Theorem}
Let $M$ be a compact quotient of $\mathbb{R} \times E(2)$. Let $g$ solve equation \eqref{bachfloweqn} where $h$ is locally homogeneous. Then $(M, g)$ converges to a flat four-dimensional manifold in the Gromov--Hausdorff topology.
\end{Theorem}

\begin{proof}
None of the components of the metric converge to zero, so there is no collapse. Since $g_{11} - g_{22} \rightarrow 0$ as $t \rightarrow \infty$, looking at the Ricci curvature, we find that the manifold becomes Ricci-flat in the limit. By Lemma \ref{sectionalcurvaturelemma}, so do the sectional curvatures.
\end{proof}

\subsection[$\widehat{{\rm SL}}(2,\mathbb{R})$]{$\boldsymbol{\widehat{{\rm SL}}(2,\mathbb{R})}$} \label{sl2rsection}

For this manifold the matrix used to determine the structure constants in Theorem~\ref{metricstructureconstantsthm} is
\[
E = \diag(-1, 1, 1).
\]
For any metric $g$, using a diagonalizing basis, the Ricci tensor is diagonal with
\begin{gather*}
{\rm Ric}_{11} = \frac{(g_{11})^2 - (g_{22} - g_{33})^2}{2 g_{22} g_{33}},\qquad
{\rm Ric}_{22} = \frac{(g_{22})^2 - (g_{11} + g_{33})^2}{2 g_{11} g_{33}},\\
{\rm Ric}_{33} = \frac{(g_{33})^2 - (g_{11} + g_{22})^2}{2 g_{11} g_{22}}
\end{gather*}
and scalar curvature is
\[
S = - \frac{(g_{11})^2 + (g_{22})^2+(g_{33})^2 + 2(g_{11} g_{22} + g_{11} g_{33} - g_{22} g_{33})}
		{2 g_{11} g_{22} g_{33}}.
\]
The Bach tensor is diagonal with
\begin{alignat*}{3}
& B_{00}= -\beta q_{\rm I}(-g_{11}, g_{22}, g_{33}) (g_{00})^3, \qquad &&
B_{11}= -\beta q_{\rm II}(-g_{11}, g_{22}, g_{33}) (g_{00})^2 g_{11},& \\
& B_{22}= -\beta q_{\rm II}(g_{22}, -g_{11}, g_{33}) (g_{00})^2 g_{22},\qquad &&
B_{33}= -\beta q_{\rm II}(g_{33}, -g_{11}, g_{22}) (g_{00})^2 g_{33},&
\end{alignat*}
where
\begin{gather*}
q_{\rm I}(x, y, z) = x^4 - x^3(y + z) + x^2yz
		 + x\big({-}y^3 + y^2z + yz^2 - z^3\big) + y^4 -y^3z -yz^3 + z^4
\end{gather*}
and
\begin{gather*}
q_{\rm II}(x, y, z) = 5x^4 - 3x^3(y + z) + x^2yz
 + x\big(y^3 - y^2z - yz^2 + z^3\big) \\
 \hphantom{q_{\rm II}(x, y, z) =}{} - 3y^4 + 3y^3z + 3yz^3 - 3z^4.
\end{gather*}
The sign choices made in the formulas for the Bach tensor here come from the fact that~$q_{\rm I}$ and~$q_{\rm II}$ are also used in the next section for $\mathbb{S}^3$, where no minus signs are needed in the expressions for the Bach tensor.

With the Bach tensor in hand, we have the following theorem:

\begin{Theorem} \label{sl2rthm}
On $\widehat{M} = \mathbb{R} \times \widehat{{\rm SL}}(2,\mathbb{R})$, every solution to equation \eqref{bachfloweqn} in a diagonalizing basis has the following properties:
\begin{itemize}\itemsep=0pt
\item $g_{00}, g_{11} \rightarrow 0$;
\item $g_{22}, g_{33} \rightarrow \infty$;
\item $g_{33} - g_{22} \rightarrow 0$.
\end{itemize}
\end{Theorem}

Before proving this, we establish some supporting lemmas. The two polynomials $q_{\rm I}$ and $q_{\rm II}$ also make an appearance in the next section so we provide some facts about them here. For $q_{\rm I}$, we have the following lemma, the proof of which is left to the reader.

\begin{Lemma} \label{ppropslemma} The polynomial $q_{\rm I}$ has the following properties:
\begin{itemize}\itemsep=0pt
\item it is symmetric;
\item $q_{\rm I}(-x, -y, -z) = q_{\rm I}(x, y, z)$;
\item it is always nonnegative;
\item it is equal to zero if and only if $x = y = z$ or two variables are equal and the third is zero.
\end{itemize}
\end{Lemma}

Note that $q_{\rm II}$ is symmetric in the last two variables. Because of this, and the fact that the flow equations for $g_{22}$ and $g_{33}$ are essentially the same, we say that without loss of generality, $g_{22} \leq g_{33}$.

The qualitative behavior of the flow is determined through a number of estimates which arise from monotonicity of various quantities. To keep things clear, these monotonicity results are presented in the following lemmas.

\begin{Lemma} \label{g33g22fraclemma}
Suppose $h_{22} < h_{33}$. Then $\frac{g_{33}}{g_{22}}$ decreases along the flow.
\end{Lemma}

\begin{proof}
Writing out the quotient rule and plugging in the differential equations for $g_{22}$ and $g_{33}$, we have
\[
\dtime \frac{g_{33}}{g_{22}}
	= -\beta \bigl(q_{\rm II}(g_{33}, -g_{11}, g_{22}) - q_{\rm II}(g_{22}, -g_{11}, g_{33}) \bigr) (g_{00})^2 \frac{g_{33}}{g_{22}}.
\]
Writing out and simplifying the factor involving the $q_{\rm II}$'s, we get
\begin{gather*}
q_{\rm II}(g_{33}, -g_{11}, g_{22}) - q_{\rm II}(g_{22}, -g_{11}, g_{33}) \\
\qquad{}= 2(g_{33}-g_{22})\big[4(g_{33})^3 + 2(g_{33})^2g_{22} + 2g_{33}(g_{22})^2 + 4(g_{22})^3 \\
 \qquad \quad {} + g_{11}\big(3(g_{33})^2 + 2g_{33}g_{22} + 3(g_{22})^2\big) + (g_{11})^3\big]
\end{gather*}
and this is positive since $g_{22} < g_{33}.$
\end{proof}

\begin{Lemma} \label{g00g22prodlemma}
Suppose $h_{22} \leq h_{33}$. Then $g_{00} g_{22}$ increases along the flow.
\end{Lemma}
\begin{proof}
Writing out the product rule and plugging in the differential equations for $g_{00}$ and $g_{22}$, we have
\[
\dtime (g_{00} g_{22})
	= -\beta \bigl(q_{\rm I}(-g_{11}, g_{22}, g_{33}) + q_{\rm II}(g_{22}, -g_{11}, g_{33})\bigr) (g_{00})^3 g_{22}.
\]
Writing out and simplifying the factor involving the $q_{\rm I}$ and $q_{\rm II}$, we get
\begin{gather*}
q_{\rm I}(-g_{11}, g_{22}, g_{33}) + q_{\rm II}(g_{22}, -g_{11}, g_{33}) \\
 \qquad{} = -2\big[(g_{33}-g_{22})\big[(g_{33})^3 + (g_{33})^2 g_{22} + g_{33} (g_{22})^2 + 3(g_{22})^3 \\
 \qquad \quad {}+ g_{11}\big((g_{33})^2 + g_{33} g_{22} + 2(g_{22})^2\big)\big]
 + (g_{11})^3 (g_{11} + g_{33})\big].
\end{gather*}
This is negative when $g_{22} \leq g_{33}$ so $\dtime (g_{00} g_{22})$ is positive.
\end{proof}

Since $g_{00}$ decreases along the flow, an immediate consequence of this lemma is that if $h_{22} \leq h_{33}$ then $g_{22}$ by itself increases along the flow. In fact, a consequence of the proof is that $\dtime g_{22}$ is positive.

\begin{Lemma} \label{g11tothepg22prodlemma} Suppose $h_{22} \leq h_{33}$. Then $(g_{11})^{\frac{3}{5}} g_{22}$ increases along the flow.
\end{Lemma}
\begin{proof}
Writing out the product rule and plugging in the differential equations for $g_{11}$ and $g_{22}$, we have
\[
\dtime \big[(g_{11})^{\frac{3}{5}} g_{22} \big]
	= -\frac{\beta}{5} \bigl(3 q_{\rm II}(-g_{11}, g_{22}, g_{33}) + 5 q_{\rm II}(g_{22}, -g_{11}, g_{33}) \bigr) (g_{00})^2 (g_{11})^{\frac{3}{5}} g_{22}.
\]
Writing out and simplifying the factor involving $q_{\rm II}$ we get
\begin{gather*}
3 q_{\rm II}(-g_{11}, g_{22}, g_{33}) + 5 q_{\rm II}(g_{22}, -g_{11}, g_{33}) \\
\qquad{} = -2\big[(g_{33}-g_{22})\big[12(g_{33})^3 + 5(g_{33})^2g_{22} + 5g_{33}(g_{22})^2 + 8(g_{22})^3 \\
\qquad \quad {}+ g_{11}(9(g_{33})^2 + 5g_{33}g_{22} + 6(g_{22})^2)\big]
 + g_{11}^2[g_{33}g_{22} + g_{11}(3g_{33} - 2g_{22})] \big].
\end{gather*}
This is always negative since $g_{22} \leq g_{33},$ so $\dtime \big[(g_{11})^{\frac{3}{5}} g_{22} \big]$ is positive.
\end{proof}

\begin{Lemma} \label{sl2difftozerolemma}Suppose $g_{22}$ and $g_{33}$ diverge and $g_{11}$ converges to zero. Then the quantity $\frac{g_{33} - g_{22}}{g_{11}}$ converges to zero along the flow.
\end{Lemma}

\begin{proof}First note that if $g_{22} = g_{33}$ then the result is true immediately. If, without loss of generality, $g_{22} < g_{33}$, we proceed in two steps. The first step is similar to the lemmas above. We have
\[
\dtime (g_{33} - g_{22}) = -\beta s(g_{11}, g_{22}, g_{33}) (g_{00})^2 (g_{33}-g_{22}),
\]
where
\begin{gather*}
s(x, y, z) = -3x^4 - x^3(y + z) - x^2yz + x\big(3y^3 + 5y^2z + 5yz^2 + 3z^3\big) \\
\hphantom{s(x, y, z) =}{} + 5y^4 + 5y^3z + 4y^2z^2 + 5yz^3 + 5z^4.
\end{gather*}
With this we have
\[
\dtime \frac{g_{33} - g_{22}}{(g_{11})^2}
	= - \beta \big[s(g_{11}, g_{22}, g_{33}) - 2q_{\rm II}(-g_{11}, g_{22}, g_{33}) \big]
		(g_{00})^2 \frac{g_{33} - g_{22}}{(g_{11})^2}.
\]
Looking at the polynomial in brackets, we have
\begin{gather*}
s(g_{11}, g_{22}, g_{33}) - 2q_{\rm II}(-g_{11}, g_{22}, g_{33}) \\
\qquad{}= g_{33}\big[(g_{33})^2 + (g_{22})^2\big](g_{33} - g_{22})
	 + 7g_{33}\big[(g_{33})^3 - (g_{11})^3\big] + 7g_{22}\big[(g_{22})^3 - (g_{11})^3\big] \\
\qquad\quad{} + 3g_{22} g_{33}\big[g_{22} g_{33} - (g_{11})^2\big]
	 + g_{11}\big[5(g_{22})^3 + 3(g_{22})^2 g_{33} + g_{22} (g_{33})^2 + 5 (g_{33})^3\big] \\
\qquad\quad{} + 3(g_{33})^4 + 4 (g_{22})^4 - 13 (g_{11})^4.
\end{gather*}
Since $g_{22}$ and $g_{33}$ go to infinity, and $g_{11}$ goes to zero, this must eventually become and stay positive and so the fraction $\frac{g_{33} - g_{22}}{(g_{11})^2}$ must eventually decrease. Since $\frac{1}{g_{11}}$ diverges, this impiles that $\frac{g_{33} - g_{22}}{g_{11}}$ must converge to zero.
\end{proof}

We are now in a position to prove Theorem~\ref{sl2rthm}. The proof requires considering a few different possibilities and ruling out any option other than what is described in the theorem.

\begin{proof}[Proof of Theorem \ref{sl2rthm}]
Without loss of generality, we may restrict our attention to flows that satisfy $g_{22} \leq g_{33}$. Now, first suppose $g_{00}$ converges to a value greater than zero, with the goal of ruling this possibility out. Consider two possibilities. Suppose first that $g_{22}$ remains bounded above. By Lemma~\ref{g33g22fraclemma}, $g_{33}$ must remain bounded as well. From this, since $\det g$ is constant, we know that $g_{11}$ remains bounded above and also below by some positive number.

By Lemma \ref{g00g22prodlemma}, $g_{22}$ is increasing, since $g_{00}$ is decreasing, so since $g_{22}$ is bounded, it must converge. By Lemma~\ref{odeconvergencelemma}, there must be a point where $\dtime g_{22} = 0$. This contradicts the fact, from the proof of Lemma~\ref{g00g22prodlemma}, that $\dtime g_{22}$ is positive in the given domain.

This implies that $g_{22}$ goes to infinity, and so must $g_{33}$. Since we are still working with the possibility that $g_{00}$ does not go to zero, we may conclude that $g_{11}$ converges to zero, again since $\det g$ is constant.

Now consider the product $g_{00} (g_{11})^{\frac{6}{5}} g_{22} g_{33}$. Note that this is equal to $(\det g)(g_{11})^{\frac{1}{5}}$ which must go to zero since $\det g$ is constant and $g_{11}$ goes to zero. On the other hand
\begin{gather*}
g_{00} (g_{11})^{\frac{6}{5}} g_{22} g_{33} = g_{00} \big[(g_{11})^{\frac{3}{5}}g_{22}\big] \big[(g_{11})^{\frac{3}{5}}g_{33}\big] \geq g_{00} \big[(g_{11})^{\frac{3}{5}}g_{22}\big]^2.
\end{gather*}
By Lemma~\ref{g11tothepg22prodlemma}, the squared factor is increasing. But this implies that $g_{00}$ must go to zero, a~contradiction.

So we may conclude that $g_{00}$ converges to zero. Knowing this, since $g_{00} g_{22}$ is increasing by Lemma~\ref{g00g22prodlemma}, $g_{22}$ and hence $g_{33}$ must both diverge to $\infty$. But then, again by Lemma~\ref{g00g22prodlemma}, $g_{00} g_{22} g_{33}$ diverges and so $g_{11}$ must go to zero since $g_{00} g_{11} g_{22} g_{33} = \det g$ is constant.

Finally, we have now established the hypotheses for Lemma \ref{sl2difftozerolemma} so we may conclude that $g_{33} - g_{22} \rightarrow 0$.
\end{proof}

With the limiting behavior of the metric established, the next step is to determine the curvature.

\begin{Proposition}
On $\widehat{M} = \mathbb{R} \times \widehat{{\rm SL}}(2, \mathbb{R})$ for every solution to equation \eqref{bachfloweqn} in a diagonalizing basis, $\frac{{\rm Ric}_{11}}{g_{11}}$ converges to~$0$, and ${\rm Ric}_{22}$ and ${\rm Ric}_{33}$ both converge to $-1$. The scalar curvature converges to~$0$ as well.
\end{Proposition}

\begin{proof}We have
\[
\frac{{\rm Ric}_{11}}{g_{11}} = \frac{(g_{11})^2 - (g_{22} - g_{33})^2}{2 g_{11} g_{22} g_{33}}.
\]
By Theorem \ref{sl2rthm}, the numerator goes to zero and, since the determinant is constant, the denominator goes to infinity.

For ${\rm Ric}_{22}$ we rewrite
\begin{gather*}
{\rm Ric}_{22} = \frac{(g_{22})^2 - (g_{11} + g_{33})^2}{2 g_{11} g_{33}}
		= \frac{g_{22} - g_{33}}{g_{11}} \frac{1}{2}\left(\frac{g_{22}}{g_{33}} + 1 \right)
			- \frac{g_{11}}{2 g_{33}}
			- 1.
\end{gather*}
By Lemma \ref{sl2difftozerolemma}, the first factor in the first term goes to zero and the rest of the term is bounded. The middle term also goes to zero. The computation for ${\rm Ric}_{33}$ is similar.

For scalar curvature, we rewrite to get
\[
S = -\frac{g_{11}}{2 g_{22} g_{33}} - \frac{g_{22}
		- g_{33}}{g_{11}} \frac{g_{22} - g_{33}}{g_{22} g_{33}}
		- \frac{1}{g_{33}} - \frac{1}{g_{22}}
\]
and, by Theorem~\ref{sl2rthm} and Lemma~\ref{sl2difftozerolemma}, all these terms go to zero.
\end{proof}

Finally, we have the following:

\begin{Theorem}Let $M$ be a compact quotient of $\mathbb{R} \times \widehat{{\rm SL}}(2, \mathbb{R})$ and let $p \in M$. Let~$g$ solve equation~\eqref{bachfloweqn} where $h$ is locally homogeneous. Then $(M, g, p)$ converges to a flat surface in the pointed Gromov--Hausdorff topology.
\end{Theorem}

\begin{proof}From the previous proposition, and by Lemma~\ref{sectionalcurvaturelemma}, the sectional curvatures all go to zero along the flow.
 \end{proof}

\subsection[$\mathbb{S}^3$]{$\boldsymbol{\mathbb{S}^3}$}

For this manifold the matrix used to determine the structure constants in Theorem~\ref{metricstructureconstantsthm} is
\[
E = {\rm id}.
\]
For any metric $g$, using a diagonalizing basis, the Ricci tensor is diagonal with
\begin{gather*}
{\rm Ric}_{11} = \frac{(g_{11})^2 - (g_{22} - g_{33})^2}{2 g_{22} g_{33}},\qquad
{\rm Ric}_{22} = \frac{(g_{22})^2 - (g_{11} - g_{33})^2}{2 g_{11} g_{33}}, \\
{\rm Ric}_{33} = \frac{(g_{33})^2 - (g_{11} - g_{22})^2}{2 g_{11} g_{22}}
\end{gather*}
and scalar curvature is
\[
S = - \frac{(g_{11})^2 + (g_{22})^2+(g_{33})^2 - 2(g_{11} g_{22} + g_{11} g_{33} + g_{22} g_{33})}
		{2 g_{11} g_{22} g_{33}}.
\]
The Bach tensor is diagonal with{\samepage
\begin{alignat*}{3}
& B_{00}= -\beta q_{\rm I}(g_{11}, g_{22}, g_{33}) (g_{00})^3 ,\qquad &&
B_{11}= -\beta q_{\rm II}(g_{11}, g_{22}, g_{33}) (g_{00})^2 g_{11}, &\\
&B_{22}= -\beta q_{\rm II}(g_{22}, g_{33}, g_{11}) (g_{00})^2 g_{22}, \qquad &&
B_{33}= -\beta q_{\rm II}(g_{33}, g_{11}, g_{22}) (g_{00})^2 g_{33},&
\end{alignat*}
where $q_{\rm I}$ and $q_{\rm II}$ were defined in Section~\ref{sl2rsection}.}

On this space, there are a variety of possibilities for Bach flow, depending on the initial conditions. To accommodate this richer structure, we break the results into a number of theorems. Because of the symmetry in the equations, we may suppose, without loss of generality, that $h_{11} \leq h_{22} \leq h_{33}$. We first analyze the cases where at least two of the initial conditions are equal. These results will begin to illustrate the complexity of the situation and begin to provide some context for the remaining cases. Ultimately, all possibilities are analyzed, culminating in Theorem~\ref{s3convergencethm}.

\begin{Theorem} \label{s3bdythm1} On $\widehat{M} = \mathbb{R} \times \mathbb{S}^3$, let $g$ solve equation~\eqref{bachfloweqn} in a diagonalizing basis with $h_{11} = h_{22} = h_{33}$. Then $g$ is static.
\end{Theorem}

Note that in this case, $\big(N^{(2)}, g^{(2)}\big)$ is a round sphere.

\begin{Theorem} \label{s3bdythm2} On $\widehat{M} = \mathbb{R} \times \mathbb{S}^3$, let $g$ solve equation~\eqref{bachfloweqn} in a diagonalizing basis with $h_{11} = h_{22} < h_{33}$ or $h_{11} < h_{22} = h_{33} < 4 h_{11}$. Then
\begin{itemize}\itemsep=0pt
\item $g_{00} \rightarrow \big(\frac{\kappa}{3} \big)^{\frac{3}{8}} (\det h)^{-\frac{1}{2}}$;
\item $g_{11}$, $g_{22}$, $g_{33} \rightarrow \big(\frac{3}{\kappa} \big)^{\frac{1}{8}} (\det h)^{\frac{1}{2}}$;
\item the components of $g_{00}$ and $g_{22}$ are related by
\[
\big(4 \det h - g_{00} (g_{22})^3\big) (g_{00})^3g_{22} = \kappa,
\]
where
\[
\kappa = \big(4 \det h - h_{00} (h_{22})^3\big) (h_{00})^3 h_{22};
\]
\item if $h_{11} = h_{22}$, then $g_{11}$ and $g_{22}$ are increasing, and $g_{33}$ is decreasing;
\item if $h_{22} = h_{33}$, then $g_{11}$ is increasing.
\end{itemize}
\end{Theorem}

\begin{Theorem} \label{s3bdythm3}
On $\widehat{M} = \mathbb{R} \times \mathbb{S}^3$, let $g$ solve equation~\eqref{bachfloweqn} in a diagonalizing basis with $4h_{11} = h_{22} = h_{33}$. Then
\begin{gather*}
g_{00}(t) = 4 (\det h) \left(\frac{1}{2^6} t + (h_{33})^2\right)^{-\frac{3}{2}}, \\
4 g_{11}(t) = g_{22}(t) = g_{33}(t) = \left(\frac{1}{2^6} t + (h_{33})^2 \right)^{\frac{1}{2}}.
\end{gather*}
\end{Theorem}

Note that in this case, $g_{00} \rightarrow 0$ and $\big(N^{(2)}, g^{(2)}\big)$ is self-similar as it expands.

\begin{Theorem} \label{s3bdythm4}
On $\widehat{M} = \mathbb{R} \times \mathbb{S}^3$, let $g$ solve equation~\eqref{bachfloweqn} in a diagonalizing basis with $4h_{11} < h_{22} = h_{33}$. Then
\begin{itemize}\itemsep=0pt
\item $g_{00}, g_{11} \rightarrow 0$;
\item $g_{22}, g_{33} \rightarrow \infty$;
\item the components of $g_{00}$ and $g_{22}$ are related by
\[
\big(g_{00} (g_{22})^3 - 4 \det h\big) (g_{00})^3 g_{22} = \kappa,
\]
where
\[
\kappa = \big(h_{00} (h_{22})^3 - 4 \det h\big) (h_{00})^3 h_{22};
\]
\item $g_{22}$ and $g_{33}$ are increasing.
\end{itemize}
\end{Theorem}

Before proving these theorems, we introduce some new structure to help with the analysis. To capitalize on the fact that $\det g$ is constant along the flow, and to exploit the symmetry among the equations for $g_{11}$, $g_{22}$, and $g_{33}$, we introduce three new variables:
\begin{gather*}
a = (g_{00})^{\frac{1}{3}} g_{11}, \qquad
b = (g_{00})^{\frac{1}{3}} g_{22}, \qquad
c = (g_{00})^{\frac{1}{3}} g_{33}
\end{gather*}
so that $abc = \det g$, and we rewrite our system using these.
We have
\begin{gather} \label{s3abcsystem}
\begin{cases}
\displaystyle \dtime g_{00} = -\beta q_{\rm I}(a, b, c) (g_{00})^{\frac{5}{3}}, \vspace{1mm}\\
\displaystyle \dtime a = -\frac{2 \beta}{3} r(a, b, c) (g_{00})^{\frac{2}{3}} a,\vspace{1mm} \\
\displaystyle \dtime b = -\frac{2 \beta}{3} r(b, a, c) (g_{00})^{\frac{2}{3}}b, \vspace{1mm}\\
\displaystyle \dtime c = -\frac{2 \beta}{3} r(c, a, b) (g_{00})^{\frac{2}{3}}c,
\end{cases}
\end{gather}
where
\begin{gather*}
r(x, y, z) = 8x^4 - 5x^3(y + z) + 2x^2yz
+ x\big(y^3 - y^2z - yz^2 + z^3\big) \\
\hphantom{r(x, y, z) =}{}
 -4y^4 + 4y^3z + 4yz^3 - 4z^4.
\end{gather*}

\begin{figure}[t!]\centering
\begin{picture}(250,250)
\put(0,0){
\includegraphics[scale = .65, clip = true, draft = false]{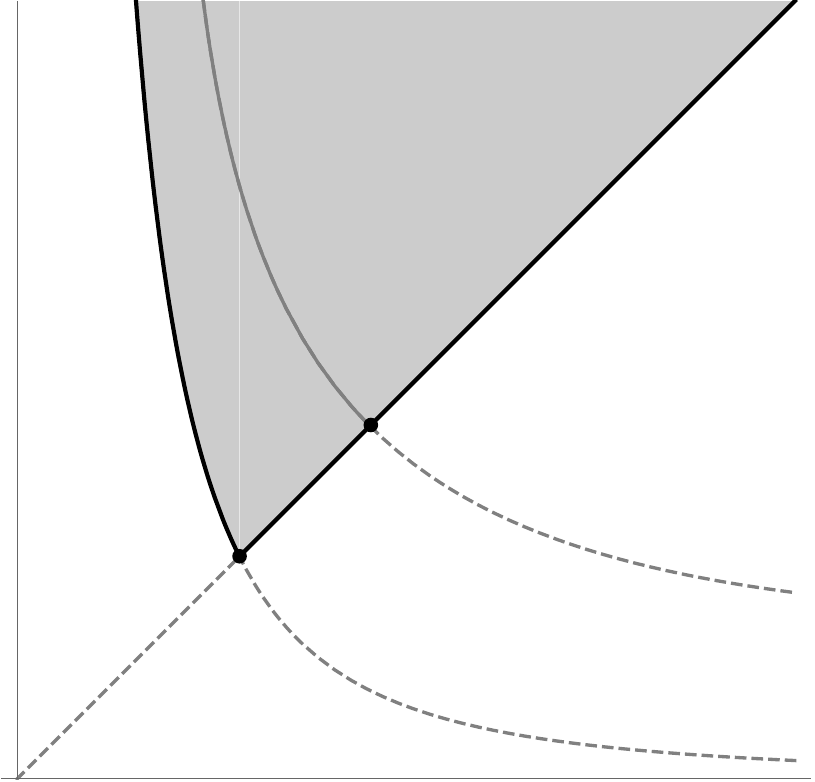}
}
\put(260,-5){$b$}
\put(5,250){$c$}
\put(72, 55){$P_0$}
\put(115, 95){$P_1$}
\put(-4, 120){$\partial D_0 \cap \partial D_{a=b}$}
\put(42, 130){\vector(2,1){14}}
\put(95, 180){$M_S$}
\put(98, 177){\vector(-1,-1){10}}
\put(170, 140){$M_U = \partial D_0 \cap \partial D_{b=c}$}
\put(178, 150){\vector(-1,2){5}}
\end{picture}
\caption{Essential features of $\partial D_0$. While the points here correspond to degenerate metrics, studying the behavior of Systems \eqref{s3abcsystem} and \eqref{s3abcadjsystem} here helps to clarify the behavior of solutions in $D$.} \label{domainfloor}
\vspace{-2mm}
\end{figure}

\begin{figure}[t!]\centering
\begin{picture}(250,250)
\put(0,0){
\includegraphics[scale = .75, viewport = 150 150 510 440, clip = true, draft = false]{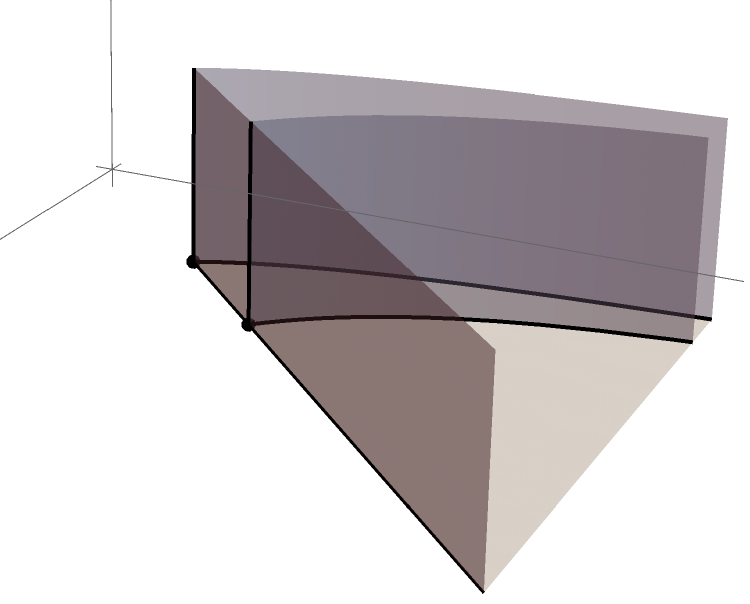}
}
\put(-5,125){$b$}
\put(275,115){$c$}
\put(35,225){$g_{00}$}
\put(55, 120){$P_0$}
\put(75, 95){$P_1$}
\put(35, 175){$L_0$}
\put(47, 180){\vector(1,0){24}}
\put(55, 150){$L_1$}
\put(67, 155){\vector(1,0){24}}
\put(240, 48){$\partial D_{0}$}
\put(241, 57){\vector(-3,2){12}}
\put(192, 40){$\partial D_{b=c}$}
\put(193, 49){\vector(-3,2){12}}
\put(150, 215){$\partial D_{a=b}$}
\put(155, 212){\vector(-1,-3){7}}
\put(200, 201){$D_S$}
\put(205, 198){\vector(-1,-3){7}}
\end{picture}
\caption{Essential features of $D$. The regions $D_{L_0}$ and $D_{\infty}$, not labelled, are determined by the surfaces shown, with $D_{L_0}$ between $\partial D_{a=b}$ and $D_S$. Together, Theorems~\ref{s3bdythm1}--\ref{s3thm3} describe the behavior of solutions in each surface or region.}\label{domain}
\vspace{-2mm}
\end{figure}

For this new system, the solution curves lie in the surface $\{a b c = \det h\}$. Moreover, because of the symmetry in the equations, we may restrict our attention to solutions that satisfy $a \leq b \leq c$.
From the determinant constraint, we know $a = \frac{\det h}{bc}$ so this inequality becomes $\sqrt{\frac{\det h}{c}} \leq b \leq c$.
Thus, the flow is analyzed on the domain
\[
D = \left\{ (g_{00}, b, c)\colon g_{00} \geq 0,\, \sqrt{\frac{\det h}{c}} \leq b \leq c \right\}.
\]
In the following, while $a$ can be eliminated, we find that it is useful to use in the analysis. As such, $a$ should always be thought of as a function of $b$ and $c$. Let $P_0$ be the point in $D$ where $g_{00} = 0$ and $a = b = c$, let $L_0$ be the ray where $g_{00} \geq 0$ and $a = b = c$, let $\partial D_0$ be the set of points in $D$ where $g_{00} = 0$, let $\partial D_{a = b}$ be the set of points in $D$ where $a = b$, let $\partial D_{b = c}$ be the set of points in $D$ where $b = c$. Note that $\partial D = \partial D_0 \cup \partial D_{a = b} \cup \partial D_{b = c}$. Let $P_1$ be the point in $D$ where $g_{00} = 0$ and $4a = b = c$, and let $L_1$ be the ray where $g_{00} \geq 0$ and $4a = b = c$. See Figs.~\ref{domainfloor} and~\ref{domain}.

With this notation in place, before proving the theorems above, we note that $\partial D_0$ corresponds to degenerate metrics and points in $\partial D_0$ are not really achievable from the perspective of the original system. However, once an initial metric $h$ is chosen, determining~$\beta$, system~\eqref{s3abcsystem} is well defined on~$\partial D_0$, and it is useful to explore the behavior here because it informs the behavior on the interior. All solutions starting here are static and as a consequence, it is conceivable that nondegenerate solutions converge to these points. We will find that with the exception of solutions converging to~$P_1$, this is not the case. In the following proofs, unless otherwise indicated, we restrict our attention to initial conditions with $g_{00} > 0$.

We now have the following:

\begin{proof}[Proof of Theorem \ref{s3bdythm1}]
This case corresponds to $L_0$. Here, $r(a, b, c) = q_{\rm I}(a, b, c) = 0$, and we have static solutions.
\end{proof}

\begin{proof}[Proof of Theorem~\ref{s3bdythm2}]
First, we consider the case where $h_{11} = h_{22}$, which corresponds to~$\partial D_{a = b}$. Here, $c = \frac{\det h}{b^2}$ and, with this, the system reduces to two variables:
\begin{gather*}
\dtime g_{00} = - \beta q_{\rm I}(b, b, c) (g_{00})^{\frac{5}{3}},\qquad
\dtime b = - \frac{2 \beta}{3} r(b, b, c) (g_{00})^{\frac{2}{3}} b.
\end{gather*}
Both $q_{\rm I}$ and $r$ simplify substantially:
\begin{gather*}
q_{\rm I}(b, b, c) = c^2 (b - c)^2,\qquad
r(b, b, c) = -c^2 (b - c) (b - 4c)
\end{gather*}
and from this, we can see that $b$ is increasing, since $b < c$. Hence $a$ is also increasing. This implies that $g_{11}$ and $g_{22}$ are also increasing. Furthermore, since $abc$ is constant, $c$ is decreasing. By Lemma \ref{odeconvergencelemma} the flow must converge to a point in either $L_0$ or $\partial D_0 \cap \partial D_{a=b}$. This second possibility will be ruled out below.

Next, since $r(b, b, c)$ is not zero, we have
\begin{equation*}
\frac{{\rm d} g_{00}}{{\rm d} b} = \frac{\dtime g_{00}}{\dtime b}
				= -\frac{3}{2} \frac{b - c}{b - 4c} \frac{g_{00}}{b}.
\end{equation*}
This is separable and we get
\begin{gather} \label{s3abequalbdyeqn}
g_{00} = \bigl[\kappa \bigl[\big(4 \det h - b^3\big) b \bigr]^{-1} \bigr]^{\frac{3}{8}},
\end{gather}
where
\begin{gather*}
\kappa = (h_{00})^{\frac{8}{3}} \big(4 \det h - (b(0))^3\big) b(0)
	 = \big(4 \det h - h_{00} (h_{22})^3\big) (h_{00})^3 h_{22}.
\end{gather*}
Substituting for $b$ in equation \eqref{s3abequalbdyeqn} and rearranging gives us the desired relationship between~$g_{00}$ and~$g_{22}$.

The relationship given by equation \eqref{s3abequalbdyeqn} shows us two things. First since
\[
g_{33} = (g_{00})^{-\frac{1}{3}} c = \det h (g_{00})^{-\frac{1}{3}}b^{-2}
\]
we can substitute and then differentiate with respect to $b$ to find that $g_{33}$ is decreasing.

Second, in the limit, we find that $g_{00}$ stays positive, so these solutions stay nondegenerate. In the limit, $a$, $b$, and $c$ converge to $(\det h)^{\frac{1}{3}}$ so
\begin{gather*}
g_{00} \rightarrow 3^{-\frac{3}{8}} \kappa (\det h)^{-\frac{1}{2}},\qquad
g_{11}, g_{22}, g_{33} \rightarrow 3^{\frac{1}{8}} \kappa^{-\frac{1}{3}} (\det h)^{\frac{1}{2}}.
\end{gather*}

Next, we consider the case where $h_{22} = h_{33}$, which, accounting for the allowable values for~$h_{11}$, corresponds to those points in $\partial D_{b = c}$ that lie between $L_0$ and $L_1$. Here, $a = \frac{\det h}{b^2}$ and, as above, the system reduces to two variables:
\begin{gather*}
\dtime g_{00} = - \beta q_{\rm I}(a, b, b) (g_{00})^{\frac{5}{3}},\qquad
\dtime b = - \frac{2 \beta}{3} r(b, a, b) (g_{00})^{\frac{2}{3}} b.
\end{gather*}
Again, both $q_{\rm I}$ and $r$ simplify substantially:
\begin{gather*}
q_{\rm I}(a, b, b) = a^2 (b - a)^2,\qquad
r(b, a, b) = -a^2 (b - a) (b - 4a).
\end{gather*}
Algebraically, this system is identical to the previous case, so the analysis is quite similar, and the resulting relationship between $g_{00}$ and $b$ is determined by the same equation \eqref{s3abequalbdyeqn}. This implies the same relationship for $g_{00}$ and $g_{22}$. Important differences arises when analyzing the qualitative behavior however. First, note that here, $\dtime b < 0$ so $b$ and $c$ are decreasing. This implies that~$a$ is increasing, and so~$g_{11}$ must be increasing as well. We cannot conclude that~$g_{22}$ and~$g_{33}$ are decreasing however, and it turns out that if the initial conditions are close enough to~$L_1$ then in fact $g_{22}$ and $g_{33}$ will increase for a while before eventually decreasing. The transition occurs when $g_{22} = g_{33} = 3g_{11}$, which is found by analyzing the equation for $\dtime g_{22}$ directly. Finally, while the qualitative behavior differs somewhat from the previous case, the limiting behavior is the same.
\end{proof}

\begin{proof}[Proof of Theorem \ref{s3bdythm3}]
This case corresponds to $L_1$. Here, $r(a, b, b) = r(b, a, b) = 0$ which implies that $\dtime a = \dtime b = \dtime c = 0$, and $g_{22} = g_{33} = 4g_{11}$ for all time. Then, using the fact that $\det g = g_{00} g_{11} g_{22} g_{33} = \frac{1}{4} g_{00} (g_{33})^3$ and focusing on the equation for $g_{33}$, we have
\begin{equation*}
\dtime g_{33} = -\beta q_{\rm II}\! \left(g_{33}, \frac{1}{4} g_{33}, g_{33} \right) (g_{00})^2 g_{33}
		= \frac{1}{2^7} (g_{33})^{-1}.
\end{equation*}
This is separable and we have
\[
g_{33} = \left(\frac{1}{2^6} t + h_{33}^2 \right)^{\frac{1}{2}}.
\]
Once this is known, the other three components are also known. We have
\begin{equation*}
g_{22} = g_{33},\qquad g_{11} = \frac{1}{4} g_{33}
\qquad \text{and}\qquad
g_{00} = 4 (\det h) \left(\frac{1}{2^6} t + h_{33}^2\right)^{-\frac{3}{2}}.\tag*{\qed}
\end{equation*}
 \renewcommand{\qed}{}
\end{proof}

\begin{proof}[Proof of Theorem \ref{s3bdythm4}]
This case corresponds to those points in $\partial D_{b = c}$ that do not lie between $L_0$ and $L_1$. Algebraically, the system is the same as for the second case in Theorem~\ref{s3bdythm2}. In this case, since $b > 4a$, $b$ and $c$ are increasing, and $a$ is decreasing, so $g_{22}$ and $g_{33}$ must be increasing as well.

Since $b$ is increasing, if it were bounded, it would have to converge to a point where \mbox{$r(b, a, b)\! =\! 0$} or where $g_{00} = 0$, by Lemma~\ref{odeconvergencelemma}. We will find below that because of the algebraic relationship between $g_{00}$ and $b$, $g_{00}$ is positive as long as $b < \infty$ so the only possibility is $r(b, a, b) = 0$. Since there are no points where this occurs other than $b = a$ and $b = 4a$, we find that $b$ and $c$, and hence $g_{22}$ and $g_{33}$ must diverge in the limit. The fact that $g_{11} \rightarrow 0$ follows from Lemma~\ref{s3azerobcinfinitylemma2} which appears later and is used for solutions starting at other points in~$D$ as well.

For the algebraic relationship between $g_{00}$ and $g_{22}$, the system is the same as for the second case in Theorem \ref{s3bdythm2}, and the analysis is essentially the same. Again, the fact that $b > 4a$ alters the formula for the trace of the solution so that instead of equation~\eqref{s3abequalbdyeqn}, we have
\begin{equation*} \label{s3bcequalbdyeqn}
g_{00} = \bigl[\kappa \bigl[\big(b^3 - 4 \det h\big) b \bigr]^{-1} \bigr]^{\frac{3}{8}},
\end{equation*}
where
\begin{gather*}
\kappa = (h_{00})^{\frac{8}{3}} ((b(0))^3 - 4 \det h) b(0)
	= (h_{00} (h_{22})^3 - 4 \det h) (h_{00})^3 h_{22}.
\end{gather*}
With this small change made, substituting for $b$ and rearranging produces the result.
\end{proof}

Our next goal is to determine the qualitative behavior of solutions with initial conditions that do not lie on the boundary. In light of the results above, we introduce a bit more notation and structure before stating the theorems. First, observe that there are no equilibria aside from those found on the boundary above. To see this, note that to have $\dtime b = \dtime c = 0$, we must have $r(b, a, c) = r(c, a, b) = 0$ and so in particular,
\[
r(c, a, b) - r(b, a, c) = 0.
\]
Writing the left side out explicitly, we have
\[
r(c, a, b) - r(b, a, c) = 3(c - b)\big(4c^3 + 2c^2b + 2cb^2 + 4b^3 - 3c^2a - 2abc - 3b^2a - a^3\big).
\]
Under the condition that $a \leq b \leq c$, we find that the large factor on the right is always positive so the only way we can have an equilibrium point is if $b = c$.

As mentioned earlier, all the points in $\partial D_{0}$ are equilibria making it difficult to determine qualitative behavior near $\partial D_{0}$. To resolve this, we adjust the system again. Specifically, we rescale the system by multiplying the right hand sides by the nonzero factor $\beta^{-1} (g_{00})^{-\frac{2}{3}}$ to produce the new system
\begin{gather} \label{s3abcadjsystem}
\begin{cases}
\displaystyle \dtime g_{00} = -q_{\rm I}(a, b, c) g_{00}, \vspace{1mm}\\
\displaystyle \dtime b = -\frac{2}{3} r(b, a, c) b, \vspace{1mm}\\
\displaystyle \dtime c = -\frac{2}{3} r(c, a, b) c.
\end{cases}
\end{gather}
The solutions to this system will just be reparameterizations of solutions to system \eqref{s3abcsystem}. Moreover, this system extends to a (mostly) nonzero system on $\partial D_0$ and, since $\dtime g_{00}$ is still zero here, solutions on this part of the boundary stay in this part of the boundary.

Restricting attention to $\partial D_0$, note that, consistent with the observations above, there are two equilibrium points $P_0$ and $P_1$. Disregarding the equation for $g_{00}$, the linearization at $P_0$ is
\[
\dtime \begin{pmatrix}
	b \\
	c
	\end{pmatrix}
		= 6 (\det h)^{\frac{4}{3}} \begin{pmatrix}
				-1 & 0 \\
				0 & -1
				\end{pmatrix}
				\begin{pmatrix}
					b \\
					c
					\end{pmatrix}
\]
and we have a stable equilibrium. The linearization at $P_1$ is				
\[
\dtime \begin{pmatrix}
	b \\
	c
	\end{pmatrix}
		= 2^{-\frac{7}{3}} 3 (\det h)^{\frac{4}{3}} \begin{pmatrix}
				-106 & 107 \\
				107 & -106
				\end{pmatrix}
				\begin{pmatrix}
					b \\
					c
					\end{pmatrix}
\]
resulting in a saddle.

The unstable manifold $M_U$ for the saddle is the line $\{b = c\}$. The stable manifold $M_S$ is a curve that approaches $P_1$ perpendicularly to $M_U$. See Fig.~\ref{domainfloor}.

Analyzing the ratio $\frac{c}{b}$, we have
\[
\dtime \frac{c}{b} = -\frac{2}{3} \bigl[r(c, a, b) - r(b, a, c) \bigr] \frac{c}{b}.
\]
As shown above, $r(c, a, b) - r(b, a, c)$ is always positive when $a < b < c$ so the fraction $\frac{c}{b}$ decreases as $t$ increases. From this, we find that the solutions approach the boundary $b = c$.

The set $\partial D_0 \backslash M_S$ comprises two components. By Lemma \ref{odeconvergencelemma}, and the fact that $\frac{c}{b}$ is decreasing, solutions starting in the component that includes $P_0$ converge to $P_0$ as $t \rightarrow \infty$ while solutions starting in the other component converge to $a = 0$, $b = c = \infty$.

Motivated by these observations, let $D_S$ be the set of points in $D$ where $g_{00} \geq 0$ and $(0, b, c) \in M_S$, and note that $D \backslash D_S$ comprises two components. Let $D_{L_0}$ be the component that includes~$L_0$, and let $D_{\infty}$ be the component that avoids~$L_0$. See Fig.~\ref{domain}.

We now have

\begin{Theorem} \label{s3thm1} On $\widehat{M} = \mathbb{R} \times \mathbb{S}^3$, let $g$ solve equation \eqref{bachfloweqn}. In a diagonalizing basis, suppose~$h$ corresponds to a point in $D_S$. Then
\begin{itemize}\itemsep=0pt
\item $g_{00} \rightarrow 0$;
\item $g_{11}, g_{22}, g_{33} \rightarrow \infty$;
\item $\frac{g_{22}}{g_{11}} \rightarrow 4$, and $\frac{g_{22}}{g_{33}} \rightarrow 1$.
\end{itemize}
\end{Theorem}

\begin{Theorem} \label{s3thm2}
On $\widehat{M} = \mathbb{R} \times \mathbb{S}^3$, let $g$ solve equation \eqref{bachfloweqn}. In a diagonalizing basis, suppose $h$ corresponds to a point in $D_{\infty}$. Then
\begin{itemize}\itemsep=0pt
\item $g_{00}, g_{11} \rightarrow 0$;
\item $g_{22}, g_{33} \rightarrow \infty$;
\item $g_{33} - g_{22} \rightarrow 0$.
\end{itemize}
\end{Theorem}

\begin{Theorem} \label{s3thm3}
On $\widehat{M} = \mathbb{R} \times \mathbb{S}^3$, let $g$ solve equation~\eqref{bachfloweqn}. In a diagonalizing basis, suppose~$h$ corresponds to a point in $D_{L_0}$. Then
\begin{itemize}\itemsep=0pt
\item $g_{00}$ does not converge to zero;
\item $g_{11}$, $g_{22}$, and $g_{33}$ converge to the same value.
\end{itemize}
\end{Theorem}

We prove these by analyzing the behavior of system~\eqref{s3abcadjsystem} and we note the following general structure for its solutions. Let $(0, b(t), c(t))$ be a solution in~$\partial D_0$ and consider the solution with initial condition $(h_{00}, b(0), c(0))$. Since the equations for $b$ and $c$ do not depend on $g_{00}$, $b(t)$ and~$c(t)$ still solve this system. Then
\[
\dtime g_{00} = -q_{\rm I}(a(t), b(t), c(t)) g_{00},
\]
which is separable and we have
\[
g_{00} = h_{00} e^{Q(t)},
\]
where
\begin{gather} \label{s3integraleqn}
Q(t) = \int_0^t -q_{\rm I}(a(\tau), b(\tau), c(\tau)) {\rm d}\tau.
\end{gather}
Note that, since $Q$ does not depend on $g_{00}$, the ratio of two solutions with initial conditions that differ only in $h_{00}$ will be constant.

We can now prove Theorems \ref{s3thm1}, \ref{s3thm2}, and \ref{s3thm3}. While Theorem \ref{s3thm1} is straightforward, it turns out that Theorems \ref{s3thm2} and \ref{s3thm3} are fairly subtle.

\begin{proof}[Proof of Theorem \ref{s3thm1}]
Since $h$ corresponds to a point in $D_S$, we know that $(0, b(0), c(0)) \in M_S$. Hence $4a(t)$, $b(t)$, and $c(t)$ all converge to the same value. The theorem will then be proved once it is established that $g_{00} \rightarrow 0$. Since the solution is bounded for $t \geq 0$, the interval on which it is defined includes $[0, \infty)$. Moreover, since the solution is converging to $L_1$, $q_{\rm I}(a(t), b(t), c(t))$ is bounded below by a positive constant, so $Q(t) \rightarrow -\infty$ as $t \rightarrow \infty$. Therefore $g_{00} \rightarrow 0$, as desired.
\end{proof}

For Theorem \ref{s3thm2}, we first establish a couple lemmas.

\begin{Lemma} \label{s3azerobcinfinitylemma1}
Suppose $a \leq b \leq c$ and $b$ and $c$ diverge $($so $a$ converges to zero$)$. Then for all $m \in \mathbb{R}$, $\frac{c - b}{a^m} \rightarrow 0$.
\end{Lemma}

\begin{proof}
Observe first that if $\frac{c - b}{a^{m'}}$ is bounded for a particular exponent $m'$, then the result is true for all $m < m'$. Therefore, it is enough to show that $\frac{c - b}{a^{m}}$ is bounded for all $m \geq 1$.

The equation solved by $c - b$ is
\[
\dtime(c - b) = -\frac{2}{3} u(a, b, c)(c - b),
\]
where
\begin{gather*}
u(x, y, z) = -4x^4 + x^3(y + z) - x^2yz -x\big(5y^3 + 7y^2 z + 7yz^2 + 5z^3\big) \\
\hphantom{u(x, y, z) =}{} + 8y^4 + 7y^3z + 6y^2z^2 + 7yz^3 + 8z^4,
\end{gather*}
and so
\[
\dtime \frac{c - b}{a^m} = -\frac{2}{3} [u(a, b, c) - m r(a, b, c)] \frac{c-b}{a^m}.
\]
Writing out $u(a, b, c) - m r(a, b, c)$, we find
\begin{gather*}
u(a, b, c) - m r(a, b, c) = 3b^3 + 6b^2 c^2 + 3c^4
					 + b^2(5b + 7c)(b - a) + c^2(7b + 5c)(c - a) \\
\hphantom{u(a, b, c) - m r(a, b, c) =}{} + m(c - b)^2\big[3b^2 + 3bc + 4c^2 + (b-a)(b+c)\big]
					 + (1+5m)a^3(b+c) \\
\hphantom{u(a, b, c) - m r(a, b, c) =}{} - a(1+2m)\big(abc + 4a^3\big).
\end{gather*}
As $a \rightarrow 0$ and $b$ and $c$ diverge to~$\infty$, the first and second lines are positive and diverge, and the third and fourth lines are positive. Only the last line is negative, but it converges to zero (since $abc = \det h$ is constant). From this, we find that $\frac{c - b}{a^m}$ is eventually decreasing, and hence bounded above.
\end{proof}

\begin{Lemma} \label{s3azerobcinfinitylemma2}
Suppose $a \leq b \leq c$ and $b$ and $c$ diverge $($so $a$ converges to zero$)$. Then $g_{11} \rightarrow 0$.
\end{Lemma}

\begin{proof}We already know that $b = (g_{00})^{\frac{1}{3}} g_{22}$ diverges. We show here that $(g_{00})^{\frac{2}{3}} g_{22}$ also eventually increases. We have
\[
\dtime \big[(g_{00})^{\frac{2}{3}} g_{22}\big]
	= -\frac{\beta}{3} s(g_{11}, g_{22}, g_{33}) (g_{00})^{\frac{8}{3}} g_{22},
\]
where
\begin{gather*}
s(x, y, z) = 2 q_{\rm I}(x, y, z) + 3 q_{\rm II}(y, x, z) \\
\hphantom{s(x, y, z)}{} = -\big[(z-y)\big(7z^3 + 6yz^2 + 6y^2z + 17y^3 - x\big(7z^2 + 6yz + 11y^2\big)\big) \\
\hphantom{s(x, y, z) =}{} + x^2\big(yz -xy - 7xz + 7x^2\big)\big].
\end{gather*}
Analyzing $s(g_{11}, g_{22}, g_{33})$, the large factor inside the first term is positive if $g_{11} \leq g_{22} \leq g_{33}$. For the second term, note that $\frac{g_{22}}{g_{11}} = \frac{b}{a}$ which diverges for the solutions under consideration. This implies that eventually, $g_{22}$ becomes, and stays, larger than $8 g_{11}$. This, combined with the fact that $g_{22} \leq g_{33}$ implies
\[
g_{22} g_{33} -g_{11} g_{22} - 7g_{11} g_{33} \geq 0.
\]
From this, we find that $s(g_{11}, g_{22}, g_{33})$ eventually becomes, and stays, negative and so $(g_{00})^{\frac{2}{3}} g_{22}$ eventually increases.

Using this fact, rewrite $\det g$ as follows:
\begin{equation*}
\det g = g_{00} g_{11} g_{22} g_{33}
	= g_{11} \big[(g_{00})^{\frac{2}{3}} g_{22}\big] c.
\end{equation*}
Since $\det g$ is constant, $(g_{00})^{\frac{2}{3}} g_{22}$ is increasing, and $c \rightarrow \infty$, it must be the case that $g_{11} \rightarrow 0$.
\end{proof}

We are now ready to prove Theorem \ref{s3thm2}

\begin{proof}[Proof of Theorem \ref{s3thm2}]
Since $h$ corresponds to a point in $D_{\infty}$, we know that $b(t)$ and $c(t)$ diverge, and so $a(t) \rightarrow 0$. Since $g_{00}$ is decreasing, it must be the case that $g_{22}$ and $g_{33}$ diverge. By Lemma~\ref{s3azerobcinfinitylemma2}, $g_{11} \rightarrow 0$, so the theorem will then be proved once it is established that $g_{00} \rightarrow 0$ as well.

For metrics in the given domain, $r(b, a, c) < 0$ so for any flow in this setting, $b$ is strictly increasing. Using this, we make a substitution to rewrite equation~\eqref{s3integraleqn} to get
\[
-Q(T) = \int_0^T q_{\rm I}(a(t), b(t), c(t)) {\rm d}t = \int_{b(0)}^{b(T)} \frac{q_{\rm I}(a, b, c)}{-\frac{2}{3} r(b, a, c) b} {\rm d} b,
\]
where we recognize that $a$ and $c$ are now functions of $b$. We now estimate $q_{\rm I}$ and $r$ along the flow. For $q_{\rm I}$, we have the following:
\begin{gather*}
q_{\rm I}(a, b, c) = a^2(b - a)(c - a) + (b - c)^2\big[a(b + c) + b^2 + bc + c^2\big]
		 \geq a^2(b - a)(c - a).
\end{gather*}
For $r$, we have
\begin{align*}
-r(b, a, c) &= a^2(b - 4a)(c - a)
		+ (c - b)\big[{-}a\big(5b^2 + 3bc + 4c^2\big) + 8b^3 + 3b^2 c + 3b c^2 + 4c^3\big] \\
		&= a^2(b - 4a)(c - a) + \frac{(c - b)}{a^5} a^5 \bar{r}(a, b, c) \\
		&= a^2\left[(b - 4a)(c - a) + \frac{(c - b)}{a^5} a^3 \bar{r}(a, b, c)\right],
\end{align*}
where $\bar{r}$ is a cubic polynomial. Since $abc$ is constant, $ab$ and $ac$ go to zero along the flow, and this implies that $a^3 \bar{r}(a, b, c) \rightarrow 0$. By Lemma \ref{s3azerobcinfinitylemma1}, $\frac{(c - b)}{a^5} \rightarrow 0$ as well. Hence, $\frac{(c - b)}{a^5} a^3 \bar{r}(a, b, c)$ is bounded by a positive constant $K$ along the flow (for $t \geq 0$) and so
\begin{gather*}
-r(b, a, c) \leq a^2([b - 4a)(c - a) + K]
		\leq L a^2(b - 4a)(c - a)
		\leq L a^2 b (c - a),
\end{gather*}
where the second inequality follows from the fact that $(b - 4a)(c - a)$ is bounded below by a~positive constant along the flow.

Combining the estimate for $q_{\rm I}$ and for $r$, we have
\begin{gather*}
\int_{b(0)}^{\infty} \frac{q_{\rm I}(a, b, c)}{-\frac{2}{3} r(b, a, c) b} {\rm d} b
		 \geq \frac{3}{2} \int_{b(0)}^{\infty} \frac{a^2(b - a)(c - a)}{L a^2 b (c - a) b} {\rm d} b
		 = \frac{3L}{2} \left( \int_{b(0)}^{\infty} \frac{1}{b} {\rm d} b - \int_{b(0)}^{\infty} \frac{a}{b^2} {\rm d} b \right).
\end{gather*}
The first integral diverges while the second integral stays finite so $Q(t) \rightarrow -\infty$ along the flow, and $g_{00} \rightarrow 0$.
\end{proof}

Before proving Theorem \ref{s3thm3}, we establish some estimates for $q_{\rm I}$ and $r$ near $P_0$.
In the following, keep in mind that since $a = \frac{\det h}{bc}$, its value changes when comparing the functions in question at different points.

\begin{Lemma} \label{s3abcequalrlemma}
There is a neighborhood $U$ of $P_0$ such that for all points $(0, b, c)$ in $U \cap \partial D_0$,
\[
r(c, a, b) \geq r(c, a, c) \geq 0.
\]
\end{Lemma}

\begin{proof}As a first step, we show that $r(c, a, b) > 0$ for points in the given domain near $a = b = c$. In fact, to help with the argument later, we show that for each $c$ the function is minimized at $b = c$. First, note that when $b = a = \sqrt{\frac{\det h}{c}}$,
\[
r(c, a, a) = 2c^2(4c - a)(c - a)
\]
and when $b = c$, so that $a = \frac{\det h}{c^2}$,
\[
r(c, a, c) = a^2(c - a)(4a - c). \]
These are both positive as long as $a < c < 4a$. Next we compute the derivative with respect to~$b$. Since $a = \frac{\deg h}{bc}$, we have $\dvar{b} a = - \frac{\det h}{c} b^{-2} = -a b^{-1}$ so
\begin{gather*}
\dvar{b} r(c, a, b) = -5c^3 - abc -3a^3 b^{-1} c - 16b^3 - 8a^3\\
\hphantom{\dvar{b} r(c, a, b) =}{} + 5ab^{-1} c^3 + 3b^2 c + 8ab^2 +16 a^4 b^{-1} + 2a^2c
\end{gather*}
and the second derivative
\begin{gather*}
\dvar{b}^2 r(c, a, b) = 12a^3 b^{-2} c - 48b^2 + 24 a^3 b^{-1} - 10ab^{-2}c^3
		 + 6bc + 8ab - 80 a^4b^{-2} - 2a^2 b^{-1} c.
\end{gather*}
This is negative at $P_0$, so must be negative in a neighborhood of this point. This implies that as $b$ varies, $r(c, a, b)$ is minimized at one of the endpoints above.

To determine which endpoint is the minimum, comparing the two expressions algebraically proves difficult. To more easily compare, let $c = v (\det h)^{\frac{1}{3}}$ (and note that $v = 1$ corresponds to the point $P_0$). Then define
\begin{gather*}
f_a(v) = r(c, a, a) = (\det h)^{\frac{4}{3}}\big(8v^4 - 10 v^{\frac{5}{2}} + 2 v\big).
\end{gather*}
Then
\[
(f_a)'(v) = (\det h)^{\frac{4}{3}}\big(32 v^3 - 25 v^{\frac{3}{2}} + 2\big)
\]
and so $(f_a)'(1) = 9 (\det h)^{\frac{4}{3}}$. Also
\[
(f_a)''(v) = (\det h)^{\frac{4}{3}} \left(96 v^2 - \frac{75}{2} v^{\frac{1}{2}} \right)
\]
and so $(f_a)''(1) = \frac{117}{2} (\det h)^{\frac{4}{3}}$.

On the other hand define
\begin{gather*}
f_c(v) = r(c, c, a) = (\det h)^{\frac{4}{3}}\big({-}4 v^{-8} + 5 v^{-5} - v^{-2}\big).
\end{gather*}
Then
\[
(f_c)'(v) = (\det h)^{\frac{4}{3}}\big(32 v^{-9} - 25 v^{-5} + 2 v^{-3}\big)
\]
and so $(f_c)'(1) = 9(\det h)^{\frac{4}{3}}$, which matches $(f_a)'(1)$. Also
\[
(f_c)''(v) = (\det h)^{\frac{4}{3}}\big({-}288 v^{-10} + 125 v^{-6} - 6 v^{-4}\big)
\]
so $(f_c)''(1) = -139 (\det h)^{\frac{4}{3}}$. This shows that $f_a$ and $f_c$ agree to first order, but that near $P_0$, $f_a$ eventually grows faster and we can conclude that for each $c$ close to $P_0$, $r(c, a, b)$ is minimized when $b = c$.
\end{proof}

\begin{Lemma} \label{s3abcequalplemma}
There is a neighborhood $U$ of $P_0$ such that for all points $(0, b, c)$ in $U \cap \partial D_0$,
\[
0 \leq q_{\rm I}(a, b, c) \leq q_{\rm I}(a, c, c).
\]
\end{Lemma}

\begin{proof}
We have
\[
q_{\rm I}(a, a, c) = c^2(c - a)^2
\qquad \text{and} \qquad
q_{\rm I}(a, c, c) = a^2(c - a)^2.
\]

For points between these two, we compute the partial derivative with respect to $b$. As in the previous lemma, $a = \frac{\deg h}{bc}$ so $\dvar{b} a = -a b^{-1}$ and we have
\begin{gather*}
\dvar{b} q_{\rm I}(a, b, c) = -4a^4 b^{-1} + 2a^3 + 3a^3 b^{-1} c - a^2 c
			 - 2ab^2 + abc + ab^{-1} c^3 + 4b^3 - 3b^2c - c^3.
\end{gather*}
The second derivative is
\begin{gather*}
\dvar{b}^2 q_{\rm I}(a, b, c) = 20a^4 b^{-2} - 6a^3 b^{-1} - 12 a^3 b^{-2} c + 2a^2 b^{-1} c
		 - 2ab - 2ab^{-2} c^3 + 12 b^2 - 6bc.
\end{gather*}
At $P_0$, this is positive so $q_{\rm I}$ is concave up near $P_0$ and we may conclude that it is maximized at one endpoint. To determine which endpoint is larger, let $c = v (\det h)^{\frac{1}{3}}$ and define
\[
f_a(v) = \bigl(q_{\rm I}(a, a, c)\bigr)^{\frac{1}{2}} = c(c - a) = (\det h)^{\frac{2}{3}} \big(v^2 - v^{\frac{1}{2}}\big)
\]
and
\[
f_c(v) = \bigl(q_{\rm I}(a, c, c)\bigr)^{\frac{1}{2}} = a(c - a) = (\det h)^{\frac{2}{3}} \big(v^{-1} - v^{-4}\big).
\]
Then $f_a(1) = f_c(1)$. Computing derivatives, we have
\[
(f_a)'(v) = (\det h)^{\frac{2}{3}} \left(2v - \frac{1}{2} v^{-\frac{1}{2}}\right)
\]
and
\[
(f_c)'(v) = (\det h)^{\frac{2}{3}} \big({-}v^{-2} + 4v^{-5}\big).
\]
Hence $(f_a)'(1) = \frac{3}{2}$ while $(f_c)'(1) = 3$ and we may conclude that near $P_0$, $f_c$ grows faster, so $q_{\rm I}$ is maximized when $b = c$.
\end{proof}

We are now ready to prove Theorem~\ref{s3thm3}.

\begin{proof}[Proof of Theorem~\ref{s3thm3}]
Since $h$ corresponds to a point in $D_{L_0}$, we know that $a(t)$, $b(t)$, and~$c(t)$ converge to the same value. The theorem will then be proved once it is established that $g_{00}$ does not go to zero.
Note that, since the solution is bounded, the interval on which it is defined includes $[0, \infty)$.

By Lemma~\ref{s3abcequalrlemma}, $r(c, b, a)$ is positive near $a = b = c$, so $\dtime c$ is negative and $c$ is strictly decreasing. From this, we can reparameterize the integral above and then use Lemmas~\ref{s3abcequalrlemma} and~\ref{s3abcequalplemma} to get
\begin{align*}
-Q(T) &= \int_0^T q_{\rm I}(a(t), b(t), c(t)) {\rm d}t
	= \int_{c(0)}^{c(T)} \frac{q_{\rm I}(a, b, c)}{-\frac{2}{3} r(c, a, b) c} {\rm d} c
	= \frac{3}{2} \int_{c(T)}^{c(0)} \frac{q_{\rm I}(a, b, c)}{r(c, a, b) c} {\rm d} c \\
	& \leq \frac{3}{2} \int_{c(T)}^{c(0)} \frac{q_{\rm I}(a, c, c)}{r(c, a, c) c} {\rm d} c
= \frac{3}{2} \int_{c(T)}^{c(0)} \frac{a^2(c-a)^2}{a^2(c-a)(4a-c) c} {\rm d} c
	= \frac{3}{2} \int_{c(T)}^{c(0)} \frac{(c-a)}{(4a-c) c} {\rm d} c.
\end{align*}
The last integrand on the second line is bounded and the interval of integration stays bounded, so the integral stays finite as $T \rightarrow \infty$. (We are allowed to cancel the factor $a^2 (c-a)$ because the~$a$ being used is the same for the numerator and the denominator, since the estimates are both taken on the same side of the boundary.)
\end{proof}

With the limiting behavior of the metric established, the next step is to determine curvature.

\begin{Proposition} \label{s3ricprop}
On $\widehat{M} = \mathbb{R} \times \mathbb{S}^3$, let $g$ be a solution to equation~\eqref{bachfloweqn} with initial metric~$h$. Then in a diagonalizing basis,
\begin{itemize}\itemsep=0pt
\item if $h$ corresponds to a point in $D_S$, then ${\rm Ric}_{11}$ converges to $\frac{1}{32}$, ${\rm Ric}_{22}$ and ${\rm Ric}_{33}$ both converge to $\frac{7}{8}$, and $S$ converges to~$0$;
\item if $h$ corresponds to a point in $D_{\infty}$, then $\frac{{\rm Ric}_{11}}{g_{11}}$ converges to zero, ${\rm Ric}_{22}$ and ${\rm Ric}_{33}$ both converge to~$1$, and $S$ converges to~$0$;
\item if $h$ corresponds to a point in $D_{L_0}$, then ${\rm Ric}_{11}$, ${\rm Ric}_{22}$, and ${\rm Ric}_{33}$ all converge to $\frac{1}{2}$, and $S$ converges to a positive value.
\end{itemize}
\end{Proposition}

\begin{proof}
Since the Ricci tensor is invariant under uniform rescaling of the metric, we can use the components of the metric directly, or we can use $a$, $b$, and $c$ to determine the Ricci curvature. We have three cases.

If $h$ corresponds to a point in $D_S$, then we can simply plug in the fact that, in the limit, $4a = b = c$ to get the values indicated for the Ricci tensor. For scalar curvature, we have{\samepage
\begin{align*}
S &= - \frac{(g_{11})^2 + (g_{22})^2+(g_{33})^2 - 2(g_{11} g_{22} + g_{11} g_{33} + g_{22} g_{33})}
			{2 g_{11} g_{22} g_{33}} \\
	&= -(g_{00})^{\frac{1}{3}} \frac{a^2 + b^2 + c^2 - 2(ab + ac + bc)}{2 \det h} \rightarrow 0
\end{align*}
since $g_{00} \rightarrow 0$.}

If $h$ corresponds to a point in $D_{\infty}$, we have
\begin{gather*}
{\rm Ric}_{33} = \frac{c^2 - (b-a)^2}{2ab}
		 = \frac{(c-b)(c+b)}{2ab} + 1 - \frac{a}{2b}.
\end{gather*}
The first term goes to zero since $\frac{c - b}{a} \rightarrow 0$ by Lemma \ref{s3azerobcinfinitylemma1}, and $\frac{c+b}{b}$ stays bounded since $\frac{c}{b} \rightarrow 1$. The last term also goes to zero so ${\rm Ric}_{33} \rightarrow 1$. The analysis for ${\rm Ric}_{22}$ is similar.

For $\frac{{\rm Ric}_{11}}{g_{11}}$ we have
\begin{equation*}
\frac{{\rm Ric}_{11}}{g_{11}} = \frac{a^2 - (c - b)^2}{2g_{11}bc}
				= \frac{g^{\frac{1}{3}}\big[a^2 - (c - b)^2\big]}{2abc}.
\end{equation*}
Here, the numerator goes to zero while the denominator stays constant.

For scalar curvature, we have
\begin{align*}
S &= - \frac{(g_{11})^2 + (g_{22})^2+(g_{33})^2 - 2(g_{11} g_{22} + g_{11} g_{33} + g_{22} g_{33})}
			{2 g_{11} g_{22} g_{33}} \\
	&= -(g_{00})^{\frac{1}{3}} \frac{a^2 + b^2 + c^2 - 2(ab + ac + bc)}{2 \det g}
	 = -(g_{00})^{\frac{1}{3}} \frac{(c - b)^2 + a^2 - 2a(c + b)}{2 \det g} \\
	&= \frac{-(g_{00})^{\frac{1}{3}}}{2 \det g}\big[(c-b)^2 + a^2\big]
		+ (g_{00})^{\frac{1}{3}} \left(\frac{1}{b} + \frac{1}{c} \right).
\end{align*}
In this form, we can see that both terms go to zero.

If $h$ corresponds to a point in $D_{L_0}$, then in the limit, $a = b = c$, and we get the desired values for the Ricci tensor. For scalar curvature,
\begin{align*}
S &= - \frac{(g_{11})^2 + (g_{22})^2+(g_{33})^2 - 2(g_{11} g_{22} + g_{11} g_{33} + g_{22} g_{33})}
			{2 g_{11} g_{22} g_{33}} \\
	&= -(g_{00})^{\frac{1}{3}} \frac{a^2 + b^2 + c^2 - 2(ab + ac + bc)}{2 \det h}
	\rightarrow \frac{3}{2} \frac{a^2 (g_{00})^{\frac{1}{3}}}{\det h}
	= \frac{3}{2 g_{11}},
\end{align*}
which is positive. 
\end{proof}

Finally, we have the following:

\begin{Theorem} \label{s3convergencethm}Let $M$ be a quotient of $\mathbb{R} \times \mathbb{S}^3$ and let $p \in M$. Let $g$ solve equation~\eqref{bachfloweqn} where~$h$ is locally homogeneous. Then
\begin{itemize}\itemsep=0pt
\item if $h$ corresponds to a point in $D_S$, then $(M, g, p)$ collapses to a flat three-dimensional manifold in the pointed Gromov--Hausdorff topology;
\item if $h$ corresponds to a point in $D_{\infty}$, then $(M, g, p)$ collapses to a flat surface in the pointed Gromov--Hausdorff topology;
\item if $h$ corresponds to a point in $D_{L_0}$, then $(M, g)$ converges to a quotient of the product of a~circle and the round sphere in the Gromov--Hausdorff topology.
\end{itemize}
\end{Theorem}

\begin{proof}These results follow from the previous proposition and Lemma~\ref{sectionalcurvaturelemma}.
\end{proof}

\subsection[$\mathbb{H}^3$]{$\boldsymbol{\mathbb{H}^3}$} This space is not a Lie group so the techniques used above do not apply. In fact, the analysis here is much simpler. There is a one parameter family of homogeneous metrics for $\mathbb{H}^3$ and they are all constant scalar multiples of the standard hyperbolic metric and hence Einstein. Therefore, by Proposition~\ref{EinsteinSliceprop}, $\mathbb{R} \times \mathbb{H}^3$ is static under Bach flow.

\subsection[$\mathbb{R} \times \mathbb{S}^2$ and $\mathbb{R} \times \mathbb{H}^2$]{$\boldsymbol{\mathbb{R} \times \mathbb{S}^2}$ and $\boldsymbol{\mathbb{R} \times \mathbb{H}^2}$}

While these spaces can be thought of as three-dimensional factors for various $1 \times 3$ products, they are more naturally viewed in terms of $2 \times 2$ products which are discussed in the next section.

\section[Bach flow on locally homogeneous $2 \times 2$ products]{Bach flow on locally homogeneous $\boldsymbol{2 \times 2}$ products} \label{22flowsec}

Bach flow on products of homogeneous surfaces was explored in \cite{DasKar}. While not new, we reproduce the analysis here for completeness because of the fact that three families of $1 \times 3$ products can also be viewed as $2 \times 2$ manifolds, namely quotients of $\mathbb{R} \times \widehat{N}$ where $\widehat{N}$ is $\mathbb{R} \times \mathbb{S}^2$, $\mathbb{R} \times\mathbb{H}^2$, or $\mathbb{R}^3$. Working as before on the universal cover, it seems at first that there are essentially six different cases to consider: $\mathbb{R}^2 \times \mathbb{R}^2$, $\mathbb{R}^2 \times \mathbb{S}^2$, $\mathbb{R}^2 \times \mathbb{H}^2$, $\mathbb{S}^2 \times \mathbb{S}^2$, $\mathbb{S}^2 \times \mathbb{H}^2$, $\mathbb{H}^2 \times \mathbb{H}^2$. However, it turns out that the Bach tensor does not distinguish between the spherical and hyperbolic slices and the analysis reduces to three cases, one of which is trivial (since the Bach tensor vanishes). To see this, because of the constancy of the scalar curvatures, equations~\eqref{Bach22eqn1} and~\eqref{Bach22eqn2} reduce to
\[
B_{\alpha \beta} = \frac{1}{24} \big(\big(S^{(1)}\big)^2 - \big(S^{(2)}\big)^2 \big) g_{\alpha \beta}
\qquad \text{and} \qquad
B_{jk} = \frac{1}{24} \big(\big(S^{(2)}\big)^2 - \big(S^{(1)}\big)^2 \big) g_{jk}.
\]
Since the scalar curvatures of the slices are squared, there is no way to distinguish between a~positively curved space and a negatively curved space.

For any of these spaces we can write a homogeneous product metric as
\[
g = f_1 g^{(1)} + f_2 g^{(2)},
\]
where $g^{(i)}$ is the standard metric for the $i$th slice and $f_i > 0$.

The constancy of the volume form along the flow implies that
\[
\gamma = f_1(t) f_2(t) = f_1(0) f_2(0)
\]
is a constant determined by the initial metric.

\subsection[$\mathbb{R}^2 \times \mathbb{R}^2$, $\mathbb{R}^2 \times \mathbb{S}^2$, and $\mathbb{R}^2 \times \mathbb{H}^2$]{$\boldsymbol{\mathbb{R}^2 \times \mathbb{R}^2}$, $\boldsymbol{\mathbb{R}^2 \times \mathbb{S}^2}$, and $\boldsymbol{\mathbb{R}^2 \times \mathbb{H}^2}$}

These spaces can also be thought of as $1 \times 3$ products and the analysis here completes the $1 \times 3$ product cases. The space $\mathbb{R}^2 \times \mathbb{R}^2$ is flat, so the metric is static under Bach flow. For the remaining two spaces, we have the following:

\begin{Theorem}Let $M$ be a compact quotient of $\mathbb{R}^2 \times \mathbb{S}^2$ or $\mathbb{R}^2 \times \mathbb{H}^2$ and let $p \in M$. Let $g$ solve equation~\eqref{bachfloweqn} where $h$ is locally homogeneous. Then $(M, g, p)$ converges to a flat surface in the pointed Gromov--Hausdorff topology.
\end{Theorem}

\begin{proof}Using the structure introduced above, the curvature of the first slice is zero and that of the second slice is
\[
S^{(2)} = \pm 2 f_2^{-1},
\]
(positive for $\mathbb{S}^2$ and negative for $\mathbb{H}^2$). Bach flow then reduces to
\begin{equation*}
\dtime f_1 = -\frac{1}{6} f_2^{-2} f_1,\qquad
\dtime f_2 = \frac{1}{6} f_2^{-1}.
\end{equation*}

The second equation is separable. Once its solution is found, it can be plugged into the first and we can solve the resulting separable equation for $f_1$. We end up with
\begin{gather*}
f_1(t) = \gamma \left(\frac{1}{3} t + f_2(0)^2 \right)^{-\frac{1}{2}},\qquad
f_2(t) = \left(\frac{1}{3} t + f_2(0)^2 \right)^{\frac{1}{2}}.
\end{gather*}

We find that the solutions are immortal but not ancient. The flat slice shrinks while the curved slice blows up and its scalar curvature goes to zero.
\end{proof}

\subsection[$\mathbb{S}^2 \times \mathbb{S}^2$, $\mathbb{S}^2 \times \mathbb{H}^2$, and $\mathbb{H}^2 \times \mathbb{H}^2$]{$\boldsymbol{\mathbb{S}^2 \times \mathbb{S}^2}$, $\boldsymbol{\mathbb{S}^2 \times \mathbb{H}^2}$, and $\boldsymbol{\mathbb{H}^2 \times \mathbb{H}^2}$}

These spaces are not $1 \times 3$ products. They are included here for completeness.

\begin{Theorem}Let $M$ be a compact quotient of $\mathbb{S}^2 \times \mathbb{S}^2$, $\mathbb{S}^2 \times \mathbb{H}^2$, or $\mathbb{H}^2 \times \mathbb{H}^2$. Let $g$ solve equation~\eqref{bachfloweqn} where $h$ is locally homogeneous. Then $(M, g)$ converges to a Bach-flat four-dimensional manifold in the Gromov--Hausdorff topology. The difference in the magnitude of the curvature of each slice converges to zero, and the scalar curvature converges to
\begin{itemize}\itemsep=0pt
\item $4 \gamma^{-\frac{1}{2}}$ on $\mathbb{S}^2 \times \mathbb{S}^2$,
\item $0$ on $\mathbb{S}^2 \times \mathbb{H}^2$,
\item $-4 \gamma^{-\frac{1}{2}}$ on $\mathbb{H}^2 \times \mathbb{H}^2$.
\end{itemize}
\end{Theorem}

\begin{proof}For these spaces the scalar curvatures are
\[
S^{(i)} = \pm 2 f_i^{-1},
\]
positive for spheres and negative for hyperbolic spaces, and the Bach tensor can be written
\begin{gather*}
B_{\alpha \beta} = \frac{1}{6} \big(f_1^{-2} - f_2^{-2}\big) f_1 \big(g^{(1)}\big)_{\alpha \beta}, \qquad
B_{jk} = \frac{1}{6} \big(f_2^{-2} - f_1^{-2}\big) f_2 \big(g^{(2)}\big)_{jk}.
\end{gather*}
Note that if $f_1(0) = f_2(0)$, then the solution is constant. Otherwise, suppose $f_1(0) < f_2(0)$. The constant $\gamma$ allows us to reduce the system~\eqref{bachfloweqn} to a single equation.
\[
\dtime f_1 = \frac{1}{6}\big(f_1^{-1} - \gamma^{-2} f_1^3\big).
\]
This is separable and we can solve to get
\begin{gather*}
f_1 = \sqrt{\gamma \tanh \left(\frac{t}{3\gamma} + \mu \right)},\qquad
f_2 = \sqrt{\gamma \coth \left(\frac{t}{3\gamma} + \mu \right)},
\end{gather*}
where
\[
\mu = \frac{1}{2} \ln \left(\frac{f_1(0) + f_2(0)}{f_2(0) - f_1(0)} \right).
\]
If $f_1(0) > f_2(0)$ the solutions are swapped. From this we see that, again, solutions are immortal, but not ancient. Here, the absolute values of the curvatures of the slices converge to the same value, and the manifold converges to a Bach-flat four-dimensional manifold. The limiting scalar curvatures can be calculated directly.
\end{proof}

\section{Comparison with Ricci flow} \label{riccompsec}

We finish with a comparison of the qualitative behavior between Ricci flow, as determined in~\cite{IJ} and~\cite{IJL}, and Bach flow. There are a number of ways this might be done. First, there is the choice of whether to include the one-dimensional component for Ricci flow. Also, there is the choice of whether or not to use volume-normalized Ricci flow or unmodified Ricci flow. For a~product metric, the Ricci tensor splits and on a one-dimensional manifold, the Ricci tensor is zero. This implies that unmodified Ricci flow on $\mathbb{S}^1 \times N$ leaves the one-dimensional component fixed and the behavior on the three-dimensional slice is the same as for Ricci flow on just $N$. For volume-normalized Ricci flow, the behavior on $\mathbb{S}^1 \times N$ will be somewhat different from that of volume-normalized Ricci flow on just $N$. First, there is a dimensional constant in the modifying term, and second, volume normalized flow on $\mathbb{S}^1 \times N$ does not preserve the volume of $N$. In the end, the differences in all these flows are somewhat cosmetic. Rescaling space and time in appropriate ways allows the solution to one of these problems to be modified so as to solve another. For a bit more detail, see the discussions in \cite{IJL} including the analysis for those cases that relate to the results in \cite{IJ}.

At first, it might seem most natural to compare Bach flow on $\mathbb{S}^1 \times N$ to volume-normalized Ricci flow on $\mathbb{S}^1 \times N$ since both flows are acting on the same space, and both flows preserve volume. However for volume-normalized Ricci flow the behavior of the one-dimensional slice depends on the scalar curvature of $N$, while by Proposition \ref{g00monotonicprop}, under Bach flow the one-dimensional slice never expands. On the other hand, for volume-normalized Ricci flow on $N$ the static solutions are Einstein. Similarly, the static solutions $g_{\mathbb{S}^1} \! + \tilde{g}$ for Bach flow are those for which $\tilde{g}$ is Einstein on $N$ by Proposition \ref{EinsteinSliceprop}. As such, the qualitative behavior can more easily be compared. We find that on most spaces, the qualitative behavior is the same, but there are two notable differences. If $N = \mathbb{S}^3$, volume-normalized Ricci flow always converges to the round sphere, while for Bach flow, the eventual qualitative behavior depends on the initial metric. If $N = \mathbb{S}^1 \times \mathbb{S}^2$, volume-normalized Ricci flow experiences curvature blow-up in finite time, while Bach flow collapses to a flat surface as $t \rightarrow \infty$.

\subsection*{Acknowledgements}

The author would like to thank Eric Bahuaud for the many valuable discussions while developing this paper, and the referees for their in-depth, candid feedback, and constructive suggestions for improvement.

\pdfbookmark[1]{References}{ref}
\LastPageEnding

\end{document}